\documentclass[preprint,11pt,3p]{elsarticle}
\usepackage{setspace}
\usepackage{graphicx}

\usepackage{amssymb}
 \usepackage{amsthm}
\usepackage[cmex10]{amsmath}
\usepackage{mathtools}
\usepackage{color}

\usepackage{lineno}

\usepackage{framed} 
\usepackage{multicol} 

\usepackage{nomencl}
\makenomenclature
\usepackage{etoolbox}
\renewcommand\nomgroup[1]{%
  \item[\bfseries
  \ifstrequal{#1}{S}{Sets}{%
  \ifstrequal{#1}{P}{Parameters}{%
  \ifstrequal{#1}{V}{Variables}{}}}%
]}

\journal{European Journal of Operational Research}

\begin{document}

\begin{frontmatter}

\title{An Integrated Market for Electricity and Natural Gas Systems with Stochastic Power Producers}

\author[label1]{Christos~Ordoudis\corref{cor1}} 
\address[label1]{Department
of Electrical Engineering, Technical University of Denmark, Kgs. Lyngby 2800, Denmark}
\address[label2]{Department of Applied Mathematics, University of M\'{a}laga, M\'{a}laga, Spain}

\cortext[cor1]{I am corresponding author}

\ead{chror@elektro.dtu.dk}

\author[label1]{Pierre~Pinson}
\ead{ppin@elektro.dtu.dk}

\author[label2]{Juan~M.~Morales}
\ead{juan.morales@uma.es}

\begin{abstract}
In energy systems with high shares of weather-driven renewable power sources, gas-fired power plants can serve as a back-up technology to ensure security of supply and provide short-term flexibility. Therefore, a tighter coordination between electricity and natural gas networks is foreseen. In this work, we examine different levels of coordination in terms of system integration and time coupling of trading floors. We propose an integrated operational model for electricity and natural gas systems under uncertain power supply by applying two-stage stochastic programming. This formulation co-optimizes day-ahead and real-time dispatch of both energy systems and aims at minimizing the total expected cost. Additionally, two deterministic models, one of an integrated energy system and one that treats the two systems independently, are presented. We utilize a formulation that considers the linepack of the natural gas system, while it results in a tractable mixed-integer linear programming (MILP) model. Our analysis demonstrates the effectiveness of the proposed model in accommodating high shares of renewables and the importance of proper natural gas system modeling in short-term operations to reveal valuable flexibility of the natural gas system. Moreover, we identify the coordination parameters between the two markets and show their impact on the system's operation and dispatch.   
\end{abstract}

\begin{keyword}
OR in energy \sep integrated energy systems \sep electricity and natural gas markets coordination \sep renewable energy \sep stochastic programming.
\end{keyword}

\end{frontmatter}




\section{Introduction}
\label{Intro}

Natural gas is considered an efficient and clean fuel that will have a key role in the future energy system. Older coal and nuclear power plants are gradually decommissioned and replaced by gas-fired power plants (GFPPs) and renewable sources of energy. The electric power sector is expected to be the main driver of natural gas consumption increase in the future \citep{U.S.EnergyInformationAdministration2016}, which will result in a tight coupling of both energy systems. In addition, the large-scale integration of uncertain and variable renewable energy production makes operational flexibility essential in energy systems. GFPPs are a well-suited flexible component for the energy system that can support other flexible sources in the power system (e.g., hydro power) but, most importantly, enhance the link with the natural gas system and the opportunity to exploit its available flexibility. In order to facilitate the operation of the future energy system, the structure of electricity and natural gas markets needs to be reconsidered. 

In the existing setup, the main inefficiencies stem from the imperfect coordination between the trading floors, as well as between the markets for those two energy commodities, i.e., electricity and natural gas. Traditionally, the electricity and natural gas markets are cleared independently and their communication is based on the interface provided by GFPPs via predefined coordination parameters, such as their fuel consumption, which depends on their dispatch in the electricity market and the price of natural gas \citep{Duenas2015}. These parameters are defined in a static way due to the imperfect coordination between the two markets. However, a static definition of these parameters does not permit the development of an efficient coordinated setup between the two markets that will allow the interaction of the two systems to flourish. In various regions around the world that highly depend on GFPPs for electricity production, several coordination schemes are considered in practice highlighting the need to further investigate this topic. For instance, ISO New England increased the coordination with the natural gas sector and studies the impact of fuel demand of GFPPs on the operation of both systems \citep{Babula2014}. A limited market coordination would result in jeopardizing available flexibility and reliability, which encourages regulators to define setups that support existing synergies. Additionally, current market designs are based on sequential independent auctions, with the day-ahead and the balancing markets as the main trading floors. This sequential approach prevents an inter-temporal coordination between these trading floors that could facilitate the integration of variable and uncertain renewables. On the other hand, a market design based on stochastic programming, as proposed in \citep{Pritchard2010} and \citep{Morales2012}, co-optimizes the day-ahead and balancing stages and is able to provide a dispatch that anticipates future balancing needs. 

Over the last years, an increased interest in studying the interaction between the electricity and natural gas systems has been raised. \cite{Geidl2007} and \cite{An2003} incorporate natural gas network constraints in the optimal power flow problem, while the unit commitment problem with natural gas security constraints is solved in \citep{Liu2009} and \citep{Li2008}. \cite{Biskas2016} use an Augmented Lagrangian method to jointly optimize the electricity and natural gas systems including unit commitment decisions. The aforementioned contributions utilize a steady-state approach to model the natural gas system which may result in suboptimal solutions when considering short-term operations. In contrast to the electricity system, natural gas can be stored in the network and moves with a lower speed than electricity. These characteristics are important as they endow the natural gas system with flexibility that can be exploited to ensure reliability; especially in cases of highly variable withdrawals by the GFPPs. \cite{Liu2011} and \cite{Chaudry2008} model the natural gas system with linepack (i.e., ability of storing natural gas in the pipelines) by solving a mixed-integer nonlinear and a nonlinear program, respectively. However, these approaches do not guarantee global optimal solutions, involve high computation times and are not suitable to be included in a market mechanism as it is hard to derive prices from them. Additionally, \cite{Correa-Posada2014} propose a linearization approach based on mixed-integer linear programming to efficiently approximate a dynamic natural gas model that is easier to solve at global optimality. \cite{Liu2009}, \cite{Li2008} and \cite{Alabdulwahab2015} consider the natural gas network constraints only to ensure feasibility and do not optimize the operational cost of the natural gas system. Hence, the optimal operation of the integrated energy system is not guaranteed. 

Several works have dealt with the presence of uncertain renewable production in coordinated electricity and natural gas systems. \cite{Alabdulwahab2015} present a stochastic programming approach, while \cite{Bai2016} utilize interval optimization. In both cases, a proper representation of the natural gas system with linepack is neglected. \cite{Qadrdan2014} and \cite{He2016} take the linepack into account in a stochastic nonlinear program and a robust co-optimization model that utilizes the alternating direction method of multipliers (ADMM) to separate electricity and natural gas systems, respectively. An alternative approach to improve the coordination between electricity and gas networks is proposed in \citep{Clegg2016}, where a flexibility metric is calculated and included in the electricity system operation to impose natural gas related constraints. Moreover, the impact of power-to-gas technology is studied in \citep{Clegg2015} along with its potential to facilitate wind power integration.   

Acknowledging the necessity for an improved coordination between electricity and natural gas short-term markets with high penetration of stochastic production, we propose a coupled clearing model that anticipates future balancing needs and optimally dispatches the integrated energy system. The proposed model is formulated as a two-stage stochastic programming problem inspired by \citep{Morales2012}. Moreover, we include linepack in the natural gas system modeling, which is taken into account both in the day-ahead and balancing stages. Comparing the proposed approach with a purely steady-state natural gas model, we highlight the significance of taking into account the linepack flexibility in real-time operation. The contributions of the paper are summarized as follows:

\begin{enumerate}
    \item Three market-clearing models are provided, ranging from the current decoupled and deterministic setups to the proposed coupled and stochastic approach. The aim is to identify and address the efficiency improvement by considering a coordinated framework between systems and trading floors. 
    \item A market design that couples the electricity and natural gas markets in the day-ahead and balancing stages, while considering the uncertainty introduced by stochastic power producers, is proposed. Moreover, an effective pricing scheme is developed and the relation of GFPPs' operating cost with the outcomes of the natural gas market is taken into account. 
    \item A tractable and linearized natural gas model with linepack is considered, while we show that this approach takes advantage of the flexibility of the natural gas network to facilitate the integration of renewables. Furthermore, we quantify and highlight the increased performance of the model in our analysis.
\end{enumerate}

The remainder of the paper is organized as follows. The coordination framework on which the market-clearing designs are based is introduced in Section \ref{CoordinationFramework}. Section \ref{NgDynamics} presents the natural gas system modeling, while the mathematical formulation of the market-clearing models is described in Section \ref{MC_models}. Section \ref{Results} demonstrates the results in a realistic case study. Finally, Section \ref{Conclusion} provides the conclusion with suggestions for future work.
\section{Coordination Framework}
\label{CoordinationFramework}

This section presents the coordination framework of electricity and natural gas markets along with the main trading floors in short-term operations. We first discuss the degrees of coordination both system-wise and in time and then formulate the market-clearing models. 

\subsection{Market designs}

\begin{figure}[ht]%
	\centering
	\includegraphics[scale=0.5]{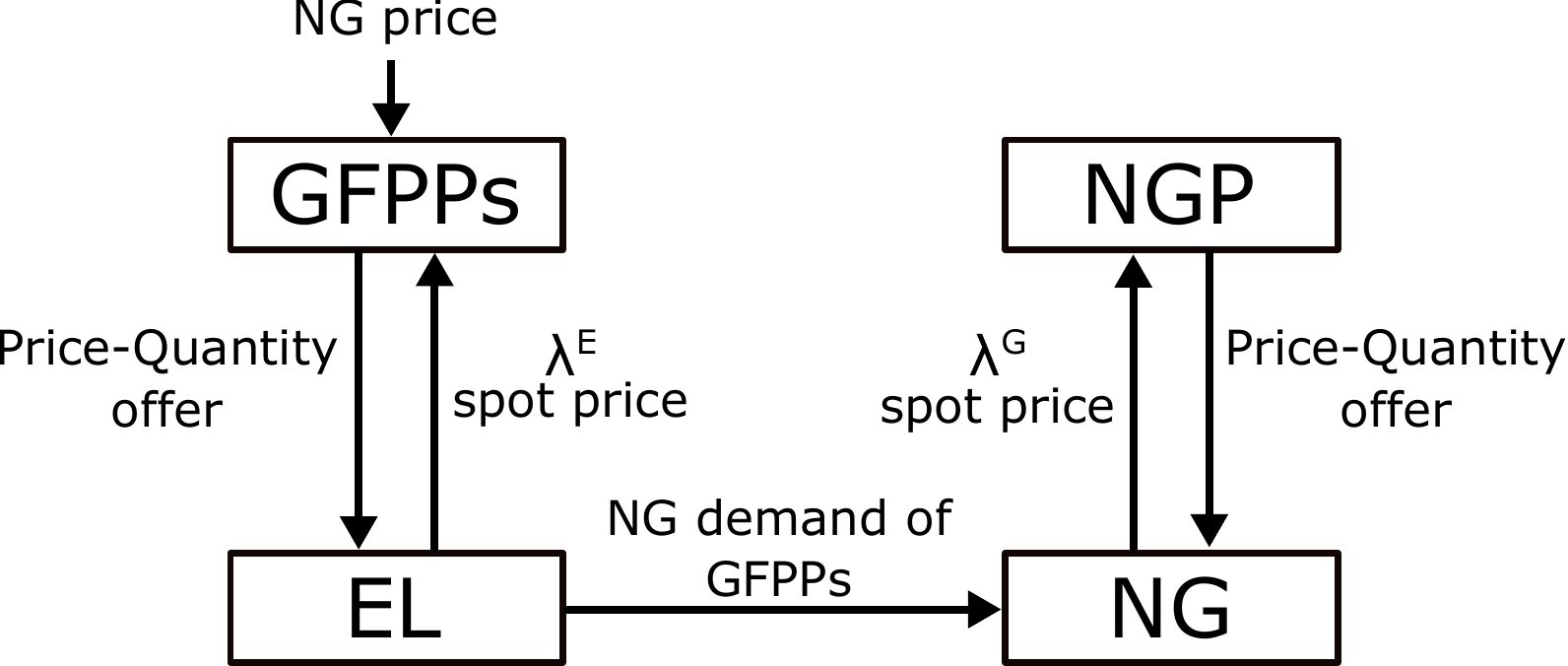}%
	\caption{Electricity and natural gas markets coordination. GFPPs: Gas-fired power plants, NGP: Natural gas producers, EL: Electricity market and NG: Natural gas market.}%
	\label{CoordinationParam}%
\end{figure}

Currently, the electricity and natural gas markets are cleared independently and mainly interact through the operation of GFPPs. This interaction is based on the definition of appropriate \textit{coordination parameters} \citep{Pardalos2012}, such as the price of natural gas and the consumption of GFPPs, as well as the available quantity of natural gas for power production. Fig. \ref{CoordinationParam} illustrates the market setup along with such coordination parameters.

Natural gas markets were historically based on long-term contracts with limited short-term variability. However, recent trends show a transition towards short-term markets, which is further strengthened by the needs of GFPPs, which are expected to have an increased intraday variability due to stochastic renewable production \citep{Zlotnik2016}. In the natural gas market, GFPPs act as buyers to acquire their fuel for power production and mainly purchase natural gas via short-term interruptible contracts or in the gas spot market. The gas supply contracts are signed with natural gas producers at a predefined price, while in the gas spot market the GFPPs buy their fuel at the spot price. The gas spot price may be undefined for a period of the following day at the time of bidding in the electricity market due to the asynchronous setup between the two markets \citep{Hibbard2012}. Consequently, GFPPs have to use an estimation of the natural gas price and face uncertainty about the price and natural gas availability, when buying from the gas spot market. The fuel price used by GFPPs to bid in the electricity market is considered fixed regardless of the procurement procedure and may not reflect the actual value of natural gas. The price-quantity offers in the electricity market are placed depending on the natural gas price, fuel availability and technical characteristics of GFPPs. Then, the electricity market is cleared and the natural gas demand of GFPPs is given as a fixed input to the natural gas market. Similarly, the natural gas market is cleared based on the price-quantity bids of natural gas producers. 

The fuel demand of GFPPs is usually assigned with a lower priority than residential natural gas loads, which may constrain the operation of GFPPs in systems where natural gas is widely used both for power production and heating, e.g., in New England \citep{Babula2014}. An additional source of inefficiency on the system operation is introduced by the asynchronous timing of electricity and natural gas markets, which becomes highly essential under an increased system interaction with fluctuating renewables. For this reason, the adoption of a concurrent market timing, as well as the consideration of the physical and economic interplay between the two markets will facilitate their coordination \citep{Hibbard2012, Tabors2012}. We capture the interaction between these markets by defining two approaches, namely the \textit{decoupled} and \textit{coupled} electricity and natural gas markets. The decoupled approach mimics the current market setup, while the coupled one addresses the need for having a concurrent and short-term integrated market. A coupled market design will profit from increased operational flexibility and provide the optimal dispatch for the whole energy system, instead of having a static and predefined communication between the two markets. A decoupled setup yields a partial coordination between markets, while a coupled approach optimizes the dispatch of GFPPs by taking into account the conditions in the natural gas market such as the price and availability of natural gas. 

In this work, the main focus is placed on two trading floors, namely the day-ahead and balancing markets, in line with the recent trend in the natural gas market where short-term trading is significantly increasing. This way, we consider a design for the natural gas market that is consistent with that of the electricity market. Current market mechanisms clear sequentially these two trading floors. The day-ahead market is initially cleared 12-36 hours in advance of actual delivery, while the balancing market settles the imbalances in relation to the day-ahead schedule to keep the system balanced \citep{Morales2014}. In such a sequential arrangement, the dispatch is based on the merit-order principle that maximizes the social welfare of each independent trading floor. However, this approach does not ensure an optimal dispatch in case more than one trading floors are evaluated. In the electricity market, flexible producers may not be scheduled or being fully dispatched, which could lead to scarcity of flexible sources in real-time operation. The aforementioned situation is aggravated with the large-scale integration of uncertain renewables. The power production cost of these sources is close to zero, which makes them the first to be dispatched according to the merit-order principle. Consequently, flexible producers may be left out of the market and higher balancing requirements will arise in real-time. Moreover, the impact of party predictable renewables on the electricity market is transferred to the natural gas market through GFPPs \citep{Keyaerts2014}. In view of an energy system where GFPPs will provide a significant amount of balancing services to support renewable production, the fuel demand of GFPPs will have a more uncertain nature. To that end, a sequential setup can be inadequate also for the natural gas market under a tighter coupling with the electricity one. 

The interaction between trading floors is modeled by a \textit{sequential} and a \textit{stochastic} approach. In the sequential setup, the day-ahead and balancing stages are optimized independently under a deterministic description of uncertainties, such as stochastic power production, which results in a lack of temporal coordination. However, a market-clearing model based on stochastic programming is able to co-optimize the dispatch in both trading floors and reduce the operating costs. The utilization of stochastic programming allows for a perfect temporal coordination between the trading floors, provided that the uncertain power production is properly described by a scenario set $\Omega$. 

The coupling of electricity and natural gas markets, along with a coordination of the two trading floors would increase the efficiency of operating the whole energy system with high shares of renewables. However, regulatory changes are still needed and market designs have to advance. The use of two-stage stochastic programming for electricity market-clearing is extensively discussed in \citep{Morales2014}, highlighting its advantages and potential challenges. On a different front, the electricity and natural gas systems have been essentially operated independently from each other over the years. It has only been recently that some regulatory changes have promoted the coordination between the two energy sectors and markets. For example, the FERC (Federal Energy Regulatory Commission) Order 809 was issued in April 16, 2015 in U.S., which aims at harmonizing the gas market with the needs of the electricity industry \citep{Zlotnik2017b}. In order to accomplish a better coordination, the information exchange needs to be improved, as well as the market structures to become more coherent. Based on this observation, we propose two coupled market-clearing models, while a decoupled one resembles current market designs that have only limited coordination. Such coordination setups have been discussed by various system operators, e.g., in New England \citep{Babula2014} and Denmark \citep{Pinson2017}, showing that it is an approach to be considered also in practice. Currently, the main inefficiencies stem from the uncoordinated operation, while most of the infrastructure already exists in both energy systems. Thus, the coupled operation of energy systems is considered inexpensive compared to solutions that involve further investments. Especially, cases like Denmark, where the two systems are operated by the same entity (i.e., Energinet.dk), would provide an environment that enables cost-efficient and coordinated operations. For these reasons, we focus on studying the three market-clearing models presented in the following section.

\subsection{Market-clearing models}

\begin{figure}[ht]%
	\centering
	\includegraphics[scale=0.5]{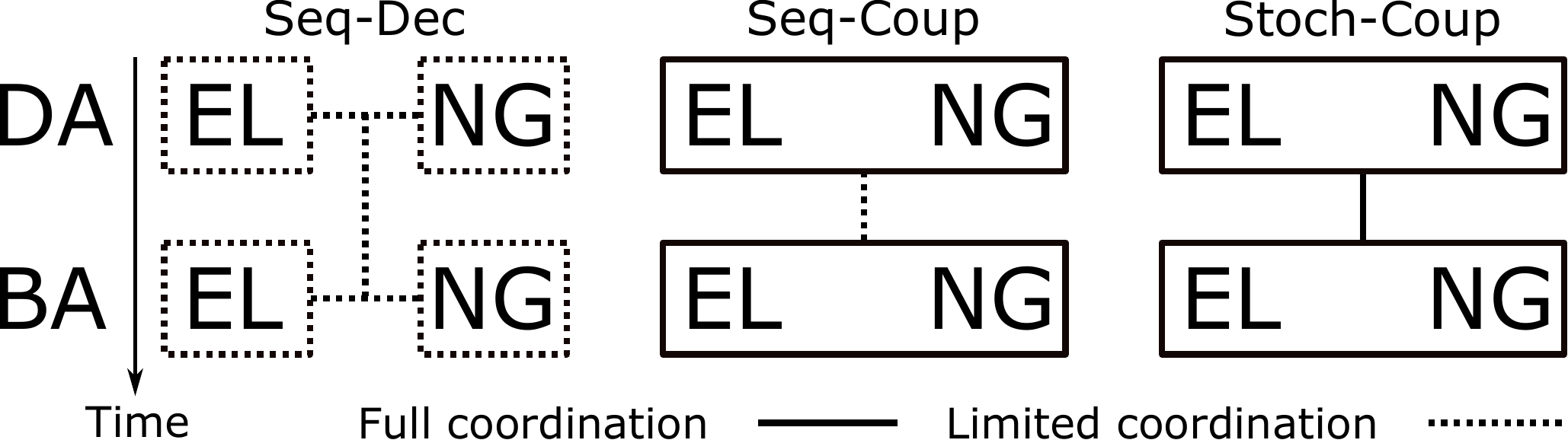}%
	\caption{Market-clearing models. DA: Day-ahead market, BA: Balancing market, EL: Electricity market and NG: Natural gas market.}%
	\label{CoordinationFig}%
\end{figure}

Depending on the level of coordination between the electricity and natural gas systems and the temporal coordination of markets, we formulate three models to cover the spectrum of potential outcomes as shown in Fig. \ref{CoordinationFig}.

\subsubsection{\textit{Model Seq-Dec -- Sequential decoupled electricity and natural gas model}} A decoupled operation of electricity and natural gas systems is considered, while renewable energy production is described in a deterministic manner. The aim of this model is to demonstrate a setup similar to the current one, where the two systems are optimized independently and the day-ahead and balancing trading floors are cleared sequentially. Initially, the electricity system is scheduled and then the fuel consumption of GFPPs is submitted as a fixed demand to the natural gas system. GFPPs have to bid in the electricity market based on an estimation of the natural gas spot price or on the price of the natural gas supply contract. However, this may not be the actual value of natural gas when operating the system. In this work, we define the \textit{actual value} of natural gas at each location of the system as its locational marginal price. Moreover, we give higher priority to residential gas loads than to the demand of GFPPs, which results in GFPPs being one of the first gas consumers to be curtailed in case it is not feasible to cover their demand. We follow an iterative approach that identifies infeasible fuel consumption schedules and constraints the allowed fuel consumption (i.e., equivalent to power production of GFPPs) for the specific time periods of the scheduling horizon \citep{Qadrdan2014}. This approach could lead to suboptimal solutions for the operation of the two systems or even infeasible schedules under cases of highly increased residential electricity and natural gas demands.

\subsubsection{\textit{Model Seq-Coup -- Sequential coupled electricity and natural gas model}} A coupled operation of electricity and natural gas systems is taken into account and renewable energy production is again described in a deterministic fashion. In this model, we formulate a single optimization problem for the operation of electricity and natural gas systems, while the day-ahead and balancing markets are still cleared sequentially. In this model, the coordination parameters are not defined in a static manner. The demand of GFPPs is treated as a variable, while the power production cost of GFPPs is defined endogenously and emerges from the cost of natural gas at the specific location. In other words, the power production cost of GFPPs is implicitly calculated through the actual value (i.e., locational marginal price) of natural gas consumed. The goal of this model is to highlight the inefficiencies that arise by having an independent scheduling of the two systems, like in \textit{Seq-Dec}.

\subsubsection{\textit{Model Stoch-Coup -- Stochastic coupled electricity and natural gas model}} The operation of electricity and natural gas systems is co-optimized and we consider a probabilistic description of uncertain renewable energy production. The proposed model is formulated as a two-stage stochastic programming problem that resembles a fully integrated system and attains perfect temporal coordination. This approach allows a pre-position of the energy system accounting for the costs of future balancing actions, in contrast to the deterministic models.

In the case that a decoupled operation of the electricity and natural gas systems is considered under uncertain renewable energy production, an additional model can be formulated. Such an approach would simulate the current decoupled setup between electricity and natural gas markets, while each operator utilizes two-stage stochastic programming to achieve a temporal coupling of trading floors. Our focus is on the integrated electricity and natural gas systems, hence the aforementioned approach is not modeled. We compare the coupled approaches with \textit{Seq-Dec}, which reproduces current market setups. The goal is to evaluate the efficiency of each model and show how the level of coordination between the electricity and natural gas systems, as well as the temporal coordination can affect the total expected cost and the dispatch of the units. 
\section{Natural Gas Network Modeling}
\label{NgDynamics}
In this section, the model of the natural gas system is presented. Natural gas system dynamics are characterized in reality by partial differential equations that are able to describe both the temporal and spatial dimensions of natural gas transport in the network \citep{Borraz-Sanchez2016}. In this work, we make some simplifying assumptions that allow us to capture some dynamics of the natural gas system, such as, the gas flow characteristics, the operation of system branches that contain compressors, the linepack flexibility and the utilization of gas storage facilities. The proposed approach takes into account the aforementioned components and allows their incorporation to the market-clearing models as day-ahead and real-time decisions. An analysis based on transient modeling of the natural gas system dynamics would provide a more realistic view of the natural gas physical behavior in the pipelines but at the same time increase the computational complexity and raise market design issues related to natural gas pricing. We refer the reader to the work of \cite{Zlotnik2017}, who incorporate transient modeling in a market framework and discuss such pricing schemes.

\subsection{Nodal and natural gas flow constraints}
The pressure at each network node has to be within specified limits in order to guarantee a secure operation of the system,
\begin{equation}
\begin{aligned}
    & PR_{m}^{\text{min}} \le pr_{m} \le PR_{m}^{\text{max}},  \enspace \forall m \in M.
 \label{Prlimits}
\end{aligned}
\end{equation}

We assume the gas transport being isothermal and in horizontal pipelines \citep{Borraz-Sanchez2016}. Moreover, the flow, pressure and linepack are defined for the middle of each pipeline as average values \citep{Correa-Posada2014}. The gas flow depends on the pressure at the adjacent nodes, the physical properties of the pipeline, such as diameter and length, as well as the volumetric characteristics of the gas. We use the Weymouth equation to describe the gas flow from node $m$ to $u$,
\begin{equation}
\begin{aligned}
    & q_{m,u} = K_{m,u}^{f}\sqrt{pr_{m}^{2}-pr_{u}^{2}},  \thinspace \forall (m,u) \in Z, \label{Weymouth}
\end{aligned}
\end{equation}
where $K_{m,u}^{f}$ is the Weymouth constant for the specific pipeline. Equation \eqref{Weymouth} is nonlinear and non-convex, thus we use an outer approximation approach based on Taylor series expansion around fixed pressure points to linearize it \citep{Tomasgard2007, Romo2009} and provide solutions that are globally optimal. Each equality constraint \eqref{Weymouth} is then replaced by a set of linear inequalities,
\begin{equation}
\begin{aligned}
    &q_{m,u}\! \le\! K_{m,u}^{f}\!\bigg(\!\frac{P\!R_{m,v} }{\sqrt{P\!R_{m,v}^{2}\!-\!P\!R_{u,v}^{2}}}pr_{\!m}\!-\!\frac{P\!R_{u,v}}{\sqrt{P\!R_{m,v}^{2}-\!P\!R_{u,v}^{2}}}pr_{\!u}\!\bigg),\\
    & \qquad \qquad \qquad \qquad \qquad \qquad \forall (m,u) \in Z, \thinspace \forall v \in V, \label{Taylor}
\end{aligned}
\end{equation}
where $V$ is the set of fixed pressure points ($P\!R_{m,v}, P\!R_{u,v}$). Using a significant number of fixed points (e.g., around 20 pairs as stated in \citep{Fodstad2015}) ensures a sufficient approximation of the Weymouth equation. The outer approximation is given by a number of planes, defined as in \eqref{Taylor}, that are tangent to the cone described by the Weymouth equation \eqref{Weymouth}. Thus, the gas flow through each of the pipelines will be approximated by the only constraint in \eqref{Taylor} that is binding \citep{Tomasgard2007}. The linear expressions used in inequalities \eqref{Taylor} are formulated based on the fixed pressure points, which in turn depend on the pressure limits of the two adjacent nodes. We generate these points by choosing multiple pressure values of the adjacent nodes between the pressure limits. The resulting fixed pressure points are used to describe the flow for one direction in the pipeline and may differ from the ones used to describe the opposite direction as the pressure limits of the adjacent nodes are not the same. For this reason, the set of inequalities \eqref{Taylor} describing the gas flow from node $m$ to $u$ may be different from the ones describing the gas flow from $u$ to $m$. We introduce two non-negative variables $q^{+}_{m,u},q^{-}_{m,u} \ge 0,\forall (m,u)\!\in\! Z$ and the binary variable $y_{m,u},\forall (m,u)\!\in\! Z$ in order to formulate the following model with bidirectional flows,
\begin{subequations}
	\begin{align}
	& q_{m,u} = q_{m,u}^{+} - q_{m,u}^{-}, \enspace  \forall (m,u) \in Z,  \label{q_linear1}\\
    & q_{m,u}^{+}\! \le\! \tilde{M}y_{m,u}, \enspace  \forall (m,u) \in Z,  \label{q_linear2}\\
    & q_{m,u}^{-}\! \le \!\tilde{M}(1-y_{m,u}), \enspace \forall (m,u) \in Z,  \enspace \label{q_linear3}\\
    & y_{m,u} + y_{u,m} = 1,  \enspace \forall (m,u) \in Z,  \label{q_linear4}\\
    & y_{m,u} \in \{0,1\},  \enspace \forall (m,u) \in Z,  \label{q_linear11}
    \end{align}
\end{subequations}
where variable $q^{+}_{m,u}$ denotes the flow in the pipeline from node $m$ to $u$, while $q^{-}_{m,u}$ from node $u$ to $m$. Parameter $\tilde{M}$ is a sufficiently large constant. Equation \eqref{q_linear1} defines the unrestricted in sign gas flow in the pipeline, while constraints \eqref{q_linear2}-\eqref{q_linear11} ensure that only one of the two variables $q^{+}_{m,u}$ and $q^{-}_{m,u}$ will take a value different from zero. In addition, we need to define the following inequalities,
\begin{subequations}
	\begin{align}
    & \big\{q_{m,u}^{+}\! \le\! KI^{+}_{m,u,v}pr_{m}-KO^{+}_{m,u,v}pr_{u} + M(1-y_{m,u}), \label{q_linear5}\\
    & q_{m,u}^{-}\! \le \!KI^{-}_{m,u,v}pr_{u}-KO^{-}_{m,u,v}pr_{m} + My_{m,u}, \label{q_linear6}\\
    & q_{u,m}^{-}\! \le\! KI^{+}_{m,u,v}pr_{m}-KO^{+}_{m,u,v}pr_{u} + My_{u,m}, \label{q_linear7}\\
    & q_{u,m}^{+}\! \le\! KI^{-}_{m,u,v}pr_{u}-KO^{-}_{m,u,v}pr_{m} + M(1-y_{u,m}) \label{q_linear8}\big\},\\
    &\nonumber \qquad \qquad \qquad \qquad \qquad \qquad \enspace \forall \{(m,u) \in Z | m < u\}, \thinspace \forall v \in V,
    \end{align}
\end{subequations}
in order to incorporate in the model the physical description of the gas flow derived by the linearization of the Weymouth equation \eqref{Taylor}. The coefficients of the linear expressions are defined as follows,

\begin{equation}
	\begin{aligned}
	& \big\{K\!I^{+}_{m,u,v}\!=\!\frac{K_{m,u}^{f}P\!R_{m,v}}{\sqrt{P\!R_{m,v}^{2}\!-\!P\!R_{u,v}^{2}}},K\!O^{+}_{m,u,v}\!=\!\frac{K_{m,u}^{f}P\!R_{u,v}}{\sqrt{P\!R_{m,v}^{2}\!-\!P\!R_{u,v}^{2}}}, \\
	& K\!I^{-}_{m,u,v}\!=\!\frac{K_{m,u}^{f}P\!R_{u,v}}{\sqrt{P\!R_{u,v}^{2}\!-\!P\!R_{m,v}^{2}}},K\!O^{-}_{m,u,v}\!=\!\frac{K_{m,u}^{f}P\!R_{m,v}}{\sqrt{P\!R_{u,v}^{2}\!-\!P\!R_{m,v}^{2}}}\big\}, \\
	& \qquad \qquad \qquad \qquad \qquad \qquad \enspace \forall \{(m,u) \in Z | m < u\}, \thinspace \forall v \in V.
	\end{aligned}
\end{equation}

It can be observed that the direction of the flow will be defined by the binary variable, which in turn will enable the appropriate linearization constraints from the set of \eqref{q_linear5}-\eqref{q_linear8}. For instance, a gas flow from node $m$ to $u$ is described by the linear inequalities that have the same coefficients, i.e., constraints \eqref{q_linear5} and \eqref{q_linear7}. It is ensured that the gas flow in the pipeline will have one direction and that $q_{m,u}\!=\!-q_{u,m},\forall (m,u)\! \in \!Z$ will hold. We also introduce two non-negative variables $q_{m,u}^{\text{in}},q_{m,u}^{\text{out}}\! \ge\! 0,\forall(m,u)\! \in \!Z$ for the inflow and outflow of each pipeline in order to account for linepack flexibility. The flow in each pipeline is defined as the average of inflow and outflow,
\begin{subequations}
	\begin{align}
    & q_{m,u}^{+}\! =\! \frac{q_{m,u}^{\text{in}}+q_{m,u}^{\text{out}}}{2}, \enspace  \forall (m,u) \in Z,  \label{q_linear9}\\
    & q_{m,u}^{-}\! =\! \frac{q_{u,m}^{\text{in}}+q_{u,m}^{\text{out}}}{2}, \enspace  \forall (m,u) \in Z.  \label{q_linear10}
    \end{align}
\end{subequations}

The branches of the network with compressors are modeled via a simplified approach that uses a compressor factor $\Gamma_{z}$ to define the relation of pressure at the two adjacent nodes as follows,
\begin{equation}
\begin{aligned}
    &  pr_{u} \le \Gamma_{z} \cdot pr_{m},  \enspace \forall (m,u) \in Z.
 \label{CompressorLimits}
\end{aligned}
\end{equation}

The inlet pressure at node $m$ can be lower than the outlet pressure at node $u$ in case the gas flows from node $m$ to $u$. The definition of gas pressures and flow through equations \eqref{Prlimits}-\eqref{CompressorLimits} is given for each pipeline $(m,u) \in Z$ of the network. However, the variables can be also indexed for each time period $t$ and scenario $\omega$ accordingly.

\subsection{Natural gas storage constraints}
An essential property of the natural gas system is that it can act as a temporary storage and be an economical way of storing energy. Linepack is very important for the short-term operation of the system and denotes the ability of storing a certain amount of natural gas in the pipeline. It is modeled by the following equations \eqref{linepack1} and \eqref{linepack2},
\begin{subequations}
	\begin{align}
	& h_{m,u,t} = K_{m,u}^{\text{h}}(pr_{m,t}+pr_{u,t})/2,  \enspace \forall (m,u) \in Z, \thinspace \forall t \in T,    \label{linepack1}\\
    & h_{m,u,t} = h_{m,u,t-1}+q_{m,u,t}^{\text{in}}-q_{m,u,t}^{\text{out}},  \enspace \forall (m,u) \in Z, \thinspace \forall t \in T. \label{linepack2}
	\end{align}
\end{subequations}

Equation \eqref{linepack1} defines the average mass of natural gas in the pipeline that is proportional to the average pressure of the adjacent nodes and a constant $K_{m,u}^{\text{h}}$ describing the pipeline characteristics. The mass conservation of each pipeline is given by \eqref{linepack2}, where the inflow and outflow of each pipeline may be different. Equations \eqref{linepack1}-\eqref{linepack2} are based on the analysis provided in \citep{Correa-Posada2014}, which describes how the pressure and the velocity of natural gas affect the mass flow. Specific levels of linepack for the initial and last optimization period should be defined to link the current day with the previous and following ones. 

Additionally, gas storage facilities are also available and can be utilized as both short- and long-term options. In short-term operation, they ensure security of supply, while they are also efficiently used as seasonal storage. The following equations model the short-term operation of gas storage facilities,
\begin{subequations}
	\begin{align}
	& E_{s}^{\text{min}} \le e_{s,t} = e_{s,t-1} + g_{s,t}^{\text{in}} - g_{s,t}^{\text{out}} \le E_{s}^{\text{max}}, \enspace \forall s \in S, \thinspace \forall t \in T, \label{storage_balance}\\
	& 0 \le g_{s,t}^{\text{in}} \le IR_{s}, \enspace \forall s \in S, \thinspace \forall t \in T, \label{storage_in}\\
    & 0 \le g_{s,t}^{\text{out}} \le WR_{s}, \enspace \forall s \in S, \thinspace \forall t \in T. \label{storage_out}
	\end{align}
\end{subequations}

Equation \eqref{storage_balance} defines the temporal balance of the gas storage facility and imposes the upper and lower limits for the volume. These limits are defined by the storage capacity and the cushion gas, which is the minimum amount of gas needed to operate the storage unit \citep{Tomasgard2007}. Constraints \eqref{storage_in}-\eqref{storage_out} enforce the maximum inflow and outflow rates.
\section{Market-clearing Models}
\label{MC_models}
The formulations of market-clearing models \textit{Seq-Dec}, $\textit{Seq-Coup}$ and \textit{Stoch-Coup} are based on the following assumptions. Regarding the stochastic power production, we take into account a single source, namely wind power, while additional uncertainties related to supply or demand side can be similarly modeled and incorporated in the models. Its stochastic nature is described through a finite set of scenarios $\Omega$ that
is available to the operator, who aims to maximize social welfare, and are properly modeled to account for temporal and spatial correlations. The uncertainty in the natural gas system stems from the random natural gas demands of GFPPs, which are the result of the stochastic wind power production. Moreover, electricity and natural gas demands are assumed to be inelastic, which allows for an equivalent formulation between the problems of social welfare maximization and cost minimization for the system operation. Thus, in this work, we take the operator's point of view who aims to minimize the total expected cost. The extreme event of load shedding is penalized with a sufficiently large value in both markets. Additionally, we assume that production costs of electricity and natural gas are described via linear functions, while the producers offer their full capacity in the market under perfect competition. More specifically, wind power producers bid in the market with zero marginal cost. The link between the two systems is provided by the GFPPs, being consumers of natural gas and producers of electricity that is sold in the electricity market. The proposed models consider network constraints in both systems. The power transmission network is modeled with DC power flow, while natural gas network modeling is also considered, as described in Section \ref{NgDynamics}. Furthermore, the models are optimized over a multi-period scheduling horizon to account for the inter-temporal constraints of the natural gas network. A concurrent market timing for the electricity and natural gas markets is assumed in the coupled models, as the energy markets are cleared simultaneously. As far as the trading floors are concerned, we consider the day-ahead and balancing stages to clear the markets, while no intraday trading is taken into account. Although this is a common assumption for the electricity market, we extend it also to natural gas market. This way the integration of electricity and natural gas system is facilitated under a consistent market design. Note that trading natural gas in short-terms markets has drawn an increased interest recently and is expected to further evolve in the future.

Following the aforementioned assumptions, the market-clearing models are recast as mixed-integer linear programming (MILP) problems. Fixing the binary variables related to the direction of the gas flow results in a linear programming (LP) problem that can properly provide day-ahead and balancing prices for electricity and natural gas as the dual variables of balancing constraints. The dual variables are indicated after a colon. A detailed mathematical description of the dispatch models is provided in the subsequent subsections.

\subsection{Sequential decoupled electricity and natural gas model}
\label{SeqDec}
The sequential and decoupled model simulates independently the operation of electricity and natural gas systems, as well as clears sequentially the day-ahead and balancing markets. Initially, the optimal day-ahead schedule that minimizes the cost of the power system is determined by model \eqref{El_DA} as follows,
\begin{subequations}
	\label{El_DA}
	\begin{align}
	& \underset{\Theta^{\text{ED}}}{\text{Min.}}\enspace \sum_{t \in T} \Big( \sum_{i_c \in I_c}C_{i_c}p_{i_c,t} + \sum_{i_g \in I_g}C_{i_g}p_{i_g,t}\Big) \label{Obj_El_DA}
	\intertext{subject to} 
	& 0 \le p_{i,t} \le P_i^{\text{max}}, \enspace \forall i \in I, \thinspace \forall t \in T, \label{El_max_min}\\
	& 0 \le w_{j,t} \le \widehat{W}_{j,t}, \enspace \forall j \in J, \thinspace \forall t \in T, \label{Wind_max_min}\\
	& \sum_{\mathclap{i \in A^{I}_{n}}} p_{i,t}\! + \!\sum_{\mathclap{j \in A^{J}_{n}}} w_{j,t}\! -\! D^{\text{E}}_{n,t} - \sum_{\mathclap{r:(n,r) \in L}}\hat{f}_{n,r,t}  = 0:\hat{\lambda}^{\text{E}}_{n,t}, \thinspace \forall n \in N, \thinspace \forall t \in T, \label{El_bal}\\
	& \hat{f}_{n,r,t} = B_{n,r}(\hat{\delta}_{n,t} - \hat{\delta}_{r,t}) \le F^{\text{max}}_{n,r}, \enspace \forall(n,r)\in L, \thinspace \forall t \in T, \label{TransCap_DA}\\
    & \hat{\delta}_{n,t} \enspace \text{free}, \enspace \forall n/n:\text{ref}, \enspace \hat{\delta}_{n,t} = 0, \enspace n: \text{ref}, \thinspace \forall t \in T, \label{RefNode_DA}
	\end{align}
\end{subequations}
where $\Theta^{\text{ED}}=\{p_{i,t},\thinspace \forall i \in I,t \in T;w_{j,t},\thinspace \forall j \in J,t \in T;\hat{\delta}_{n,t},\thinspace \forall n \in N, \thinspace \forall t \in T\}$ is the set of primal optimization variables. The objective function \eqref{Obj_El_DA} to be minimized is the operating cost of the power system, which originates from the energy production cost of the power plants. The GFPPs use an estimation of the natural gas spot price or the price of natural gas contract in order to calculate their price-quantity offers. Constraints \eqref{El_max_min} set the bounds of power production for dispatchable power plants, while constraints \eqref{Wind_max_min} limit wind power by the expected wind power production. The power balance at each node of the power system is enforced by \eqref{El_bal}. The transmission capacity limits are imposed by \eqref{TransCap_DA} and equation \eqref{RefNode_DA} defines the voltage angle at each node of the system. Having determined the day-ahead dispatch of the electricity system, we calculate the fuel consumption of GFPPs which in turn translates into a time-varying demand for each node of the natural gas system through $D_{m,t}^{\text{P}}=\sum_{i_g \in A^{I_g}_{m}} \phi_{i_g}p_{i_g,t}, \forall m \in M, \thinspace \forall t \in T$. Demand $D_{m,t}^{\text{P}}$ is treated as a parameter in the following model \eqref{Gas_DA} that defines the optimal day-ahead schedule of the natural gas system as follows,
\begin{subequations}
	\label{Gas_DA}
	\begin{align}
	& \underset{\Theta^{\text{GD}}}{\text{Min.}}\enspace \sum_{t \in T} \Big( \sum_{k \in K}C_{k}g_{k,t} + \sum_{s \in S}C_{s}g_{s,t}^{\text{out}}\Big) \label{Obj_Gas_DA}
	\intertext{subject to} 
	& 0 \le g_{k,t} \le G_{k}^{\text{max}}, \enspace \forall k \in K, \thinspace \forall t \in T, \label{Gas_max_min}\\
	& \sum_{\mathclap{k \in A^{K}_{m}}}g_{k,t} + \sum_{\mathclap{s \in A_{m}^{S}}}( g_{s,t}^{\text{out}} - g_{s,t}^{\text{in}}) - D^{\text{G}}_{m,t} - D_{m,t}^{\text{P}} \nonumber \\ 
	&- \sum_{\mathclap{u:(m,u) \in Z}}( q_{m,u,t}^{\text{in}}-q_{u,m,t}^{\text{out}})= 0:\hat{\lambda}^{\text{G}}_{m,t}, \enspace \forall m \in M, \thinspace \forall t \in T, \label{Gas_bal}	\\
	& \text{NG flow constraints} \enspace \eqref{Prlimits},\! \eqref{q_linear1}\!-\!\eqref{q_linear11},\! \eqref{q_linear5}\!-\!\eqref{q_linear8}, \eqref{q_linear9},\!\eqref{q_linear10},\! \eqref{CompressorLimits},\thinspace \forall t \in T, \label{DA_constraits_set1} \\
	& \text{NG storage constraints} \enspace \eqref{linepack1},\eqref{linepack2}, \eqref{storage_balance}-\eqref{storage_out}, \label{DA_constraits_set2}
	\end{align}
\end{subequations}
where $\Theta^{\text{GD}}=\{g_{k,t},\thinspace \forall k \in K, \thinspace \forall t \in T;e_{s,t},g_{s,t}^{\text{in}},g_{s,t}^{\text{out}}, \thinspace \forall s \in S, \thinspace \forall t \in T;h_{m,u,t},q_{m,u,t}^{\text{in}},$ $q_{m,u,t}^{\text{out}},y_{m,u,t},$ $\thinspace \forall (m,u) \in Z,$ $ \thinspace \forall t \in T;$ $pr_{m,t}, \thinspace \forall m \in M, \thinspace \forall t \in T \}$ is the set of primal optimization variables. The aim is to minimize the operating cost of the natural gas system that is represented in the objective function \eqref{Obj_Gas_DA} as the cost of natural gas production and withdrawal cost from the storage facilities. Constraints \eqref{Gas_max_min} enforce the limits of natural gas production. The gas balance for each node of the system is given by \eqref{Gas_bal}. The operation of the natural gas system is described by the set of constraints \eqref{DA_constraits_set1} and \eqref{DA_constraits_set2}. In real-time operation, wind power production $W_{j,\omega',t}$ is realized and the balancing markets of electricity and natural gas are cleared independently. The day-ahead schedule of electricity and natural gas systems is treated as fixed parameters (denoted with superscript `*') in the following formulations. The electricity balancing market \eqref{El_RT} writes as follows,
\begin{subequations}
\label{El_RT}
	\begin{align}
	& \underset{\Theta^{\text{ER}}}{\text{Min.}}\enspace \sum_{t \in T} \Big( \sum_{i_c \in I_c}(C_{i_c}^{+}p_{i_c,\omega',t}^{+}\!-\!C_{i_c}^{-}p_{i_c,\omega',t}^{-})\! +\! \sum_{n \in N}C^{\text{sh,E}}l^{\text{sh,E}}_{n,\omega',t} \nonumber\\
	& + \sum_{i_g \in I_g}(C_{i_g}^{+}p_{i_g,\omega',t}^{+}-C_{i_g}^{-}p_{i_g,\omega',t}^{-})
    \Big) \label{Obj_El_RT} 
\intertext{subject to}
& \Delta p_{i,\omega',t} = p_{i,\omega',t}^{+} - p_{i,\omega',t}^{-}, \enspace \forall i \in I, \thinspace \forall t \in T, \label{RT_def}\\
& -p^{\text{*}}_{i,t} \le \Delta p_{i,\omega',t} \le P_i^{\text{max}} -p^{\text{*}}_{i,t}, \enspace \forall i \in I, \thinspace \forall t \in T, \label{RT_capacity1}\\
& -P_i^{\text{-}} \le \Delta p_{i,\omega',t} \le P_i^{\text{+}}, \enspace \forall i \in I, \thinspace \forall t \in T, \label{RT_capacity2}
\end{align}
\begin{align}
& 0 \le l^{\text{sh,E}}_{n,\omega',t} \le D^{\text{E}}_{n,t} , \enspace \forall n \in N, \thinspace \forall t \in T, \label{RT_Lsh}\\
& 0 \le w^{\text{sp}}_{j,\omega',t} \le W_{j,\omega',t}, \enspace \forall j \in J, \thinspace \forall t \in T, \label{RT_Wspill}\\
& \sum_{\mathclap{i \in A^{I}_{n}}} \Delta p_{i,\omega',t} + l^{\text{sh,E}}_{n,\omega',t} +  \sum_{\mathclap{j \in J}}(W_{j,\omega',t} - w^{\text{sp}}_{j,\omega',t} - w^{\text{*}}_{j,t})\nonumber \\
& + \sum_{\mathclap{r:(n,r) \in L}}\hat{f}^{\text{*}}_{n,r,t}-\tilde{f}_{n,r,\omega',t} = 0:\tilde{\lambda}^{\text{E}}_{n,\omega',t}, \enspace \forall n \in N, \thinspace \forall t \in T,  \label{RT_El_balance}\\
& \tilde{f}_{n,r,\omega',t} \!=\! B_{n,r}(\tilde{\delta}_{n,\omega',t} \!-\! \tilde{\delta}_{r,\omega',t}) \!\le\! F^{\text{max}}_{n,r}, \thinspace \forall (n,r)\in L, \thinspace \forall t \in T, \label{TransCap_RT}\\
& \tilde{\delta}_{n,\omega',t} \enspace \text{free}, \enspace \forall n/n:\text{ref}, \enspace \tilde{\delta}_{n,\omega',t} = 0, \enspace n: \text{ref},\enspace \forall t \in T, \label{RefNode_RT}
\end{align}
\end{subequations}
where $\Theta^{\text{ER}}=\{p_{i,\omega',t}^{+}, p_{i,\omega',t}^{-},\thinspace \forall i \in I, \thinspace \forall t \in T; w^{\text{sp}}_{j,\omega',t},\thinspace \forall j \in J, \thinspace \forall t \in T;\tilde{\delta}_{n,\omega',t},l^{\text{sh,E}}_{n,\omega',t}, \thinspace \forall n \in N, \thinspace \forall t \in T\}$ is the set of primal optimization variables. The objective function \eqref{Obj_El_RT} describes the cost of re-dispatch actions, i.e., up/down regulation and load shedding. The aim is to activate the least-cost re-dispatch actions in order to maintain the system in balance. Power regulation is defined by equation \eqref{RT_def} and constraints \eqref{RT_capacity1} determine its bounds by taking into account the day-ahead dispatch of power plants. Moreover, constraints \eqref{RT_capacity2} limit power regulation by the reserve capacity offer. These reserve capacity offers are defined via reserve capacity markets and are treated as parameters in the market-clearing models presented in this work. Reserve capacity markets are cleared independently but can be incorporated in the presented market-clearing models as described in \citep{Morales2012}. Electricity load shedding and wind spillage are constrained by the nodal demand and actual wind power realization though \eqref{RT_Lsh} and \eqref{RT_Wspill}, respectively. Constraints \eqref{RT_El_balance} guarantee power balance at each node of the electricity network in real-time operation. Constraints \eqref{TransCap_RT} enforce power transmission capacity limits, while equation \eqref{RefNode_RT} defines the voltage angle of the system nodes. Similarly to the day-ahead stage, the fuel consumption of the GFPPs is converted to a time-varying demand deviation via $D_{m,\omega',t}^{\text{PR}}=\sum_{i_g \in A^{I_g}_{m}} \phi_{i_g}\Delta p_{i_g,\omega',t}, \forall m \in M, \thinspace \forall t  \in T$. The balancing natural gas market is formulated in \eqref{Gas_RT} as follows,
\begin{subequations}
\label{Gas_RT}
	\begin{align}
	& \underset{\Theta^{\text{GR}}}{\text{Min.}}\enspace \sum_{t \in T} \Big( \sum_{k \in K}(C_{k}^{+}g_{k,\omega',t}^{+}-C_{k}^{-}g_{k,\omega',t}^{-}) + \sum_{m \in M} C^{\text{sh,G}}l^{\text{sh,G}}_{m,\omega',t} \nonumber\\ 
	& +\sum_{s \in S}(C_{s}^{+}g_{s,\omega',t}^{\text{out},r}-C_{s}^{-}g_{s,\omega',t}^{\text{in},r}) \Big) \label{Obj_Gas_RT}
\intertext{subject to}
& 0 \le g_{k,\omega',t}^{+} \le G_{k}^{\text{max}}-g^{\text{*}}_{k,t}, \enspace \forall k \in K, \thinspace \forall t \in T,\label{RT_capacityG1}\\
& 0 \le g_{k,\omega',t}^{-} \le g^{\text{*}}_{k,t}, \enspace \forall k \in K, \thinspace \forall t \in T,\label{RT_capacityG2}\\
& 0 \le g_{k,\omega',t}^{+} \le G^{+}_{k}, \enspace \forall k \in K, \thinspace \forall t \in T,\label{RT_RcapacityG1}\\
& 0 \le g_{k,\omega',t}^{-} \le G^{-}_{k}, \enspace \forall k \in K, \thinspace \forall t \in T,\label{RT_RcapacityG2}
\end{align}
\begin{align}
& 0 \le l^{\text{sh,G}}_{m,\omega',t} \le D^{\text{G}}_{m,t}, \enspace \forall m \in M, \thinspace \forall t \in T, \label{Gas_shed}\\
& \sum_{\mathclap{k \in A^{K}_{m}}} (g_{k,\omega',t}^{+} - g_{k,\omega',t}^{-})+ \sum_{\mathclap{s \in A_{m}^{S}}}( g_{s,\omega',t}^{\text{out,r}} - g_{s,\omega',t}^{\text{in,r}})   + \sum_{\mathclap{u:(m,u) \in Z}}(q_{m,u,t}^{\text{in,*}}-q_{u,m,t}^{\text{out,*}} - q_{m,u,t}^{\text{in,r}}+q_{u,m,t}^{\text{out,r}}) \nonumber \\
& - D_{m,\omega',t}^{\text{PR}}  + l^{\text{sh,G}}_{m,\omega',t} = 0:\tilde{\lambda}^{\text{G}}_{m,\omega',t}, \enspace \forall m \in M, \thinspace \forall t \in T,  \label{RT_Gas_balance}\\
& \text{NG flow constraints} \enspace \eqref{Prlimits},\! \eqref{q_linear1}\!-\!\eqref{q_linear11},\! \eqref{q_linear5}\!-\!\eqref{q_linear8}, \eqref{q_linear9},\!\eqref{q_linear10},\! \eqref{CompressorLimits},\thinspace \forall t \in T, \label{RT_constraits_set1}\\
& \text{NG storage constraints} \enspace \eqref{linepack1},\eqref{linepack2}, \eqref{storage_balance}-\eqref{storage_out}, \label{RT_constraits_set2}
\end{align}
\end{subequations}
where $\Theta^{\text{GR}}=\{g_{k,\omega',t}^{+}, g_{k,\omega',t}^{-},\thinspace \forall k \in K, \thinspace \forall t \in T;e^{\text{r}}_{s,\omega',t},g_{s,\omega',t}^{\text{in,r}},$ $g_{s,\omega',t}^{\text{out,r}}, \thinspace \forall s \in S, \thinspace \forall t \in T;$ $h^{\text{r}}_{m,u,\omega',t},$ $q_{m,u,\omega',t}^{\text{in,r}},q_{m,u,\omega',t}^{\text{out,r}},$ $y^{\text{r}}_{m,u,\omega',t}\thinspace \forall (m,u) \in Z, \thinspace \forall t \in T;$ $pr^{\text{r}}_{m,\omega',t},l^{\text{sh,G}}_{m,\omega',t}, \thinspace \forall m \in M, \thinspace \forall t \in T\}$ is the set of primal optimization variables. In order to maintain the natural gas system balanced, a set of re-dispatch actions can be activated, namely, up/down regulation from the producers, regulation from the storage facilities and load shedding. The costs of these actions comprise objective function \eqref{Obj_Gas_RT}. The upward and downward regulation from natural gas producers is limited by constraints \eqref{RT_capacityG1} and \eqref{RT_capacityG2} that incorporate also the day-ahead schedules. Additionally, constraints \eqref{RT_RcapacityG1} and \eqref{RT_RcapacityG2} ensure that the regulation from natural gas producers is kept above zero and below the reserve capacity offer. Similarly to the electricity market, reserve capacity offers are defined through the corresponding reserve capacity markets. Natural gas load shedding is limited by the nodal demand through \eqref{Gas_shed}. Observe that an unexpected deviation of wind power production from the day-ahead schedule has to be covered by a re-dispatch action in the power system, which in turn may translate into a deviation for the natural gas system through parameter $D_{m,\omega',t}^{\text{PR}}$ that is also indexed by scenario $\omega'$. Constraint \eqref{RT_Gas_balance} guarantees the gas balance at each node of the system in real-time operation. The first re-dispatch action that is practically free stems from the linepack flexibility considered in this model through the representation of the natural gas system via constraints \eqref{RT_constraits_set1} and \eqref{RT_constraits_set2}. This is also a common practice in reality as the natural gas system is mainly balanced via linepack and the costly re-dispatch actions are only activated when flexibility by the temporal coupling of linepack is not available. It should be noted that variables in constraints \eqref{RT_constraits_set1} and \eqref{RT_constraits_set2} are to be augmented with superscript $r$ and indexed by scenario $\omega'$. We follow an iterative approach to simulate the operation of electricity and natural gas system in the case of extreme events, such as gas load shedding. If shedding of residential natural gas load takes place, we identify the GFPPs that cause this event and the corresponding nodes of the system. Then, the balancing electricity market is cleared again with additional constraints that enforce bounds on the maximum power production of GFPPs during the necessary time periods. This procedure is repeated until no residential gas load shedding occurs due to the fuel consumption of GFPPs. Solving models \eqref{El_DA}-\eqref{Gas_RT} results in calculating the total expected cost of the system under scenario set $\Omega$.
\subsection{Sequential coupled electricity and natural gas model}
\label{SeqCoup}
The sequential and coupled dispatch model simulates an energy system where the day-ahead and balancing markets are sequentially cleared, while the operation of electricity and natural gas systems is coordinated. The optimal schedule that minimizes the day-ahead cost of the integrated system is determined by model \eqref{El_Gas_DA} as follows,
\begin{subequations}
	\label{El_Gas_DA}
	\begin{align}
	& \underset{\Theta^{\text{D}}}{\text{Min.}}\enspace \sum_{t \in T} \Big( \sum_{i_c \in I_c}C_{i_c}p_{i_c,t} + \sum_{k \in K}C_{k}g_{k,t} + \sum_{s \in S}C_{s}g_{s,t}^{\text{out}}\Big) \label{Obj_DA_coup}
	\intertext{subject to} 
    & \text{DA EL constraints} \enspace \eqref{El_max_min}-\eqref{RefNode_DA}, \label{DA_constraints_coup1}\\
    & \text{DA NG constraints} \enspace \eqref{Gas_max_min},\eqref{DA_constraits_set1},\eqref{DA_constraits_set2}, \label{DA_constraints_coup2}\\
    & \sum_{\mathclap{k \in A^{K}_{m}}}g_{k,t} + \sum_{\mathclap{s \in A_{m}^{S}}}( g_{s,t}^{\text{out}} - g_{s,t}^{\text{in}}) - D^{\text{G}}_{m,t} - \sum_{\mathclap{i_g \in A^{I_g}_{m}}} \phi_{i_g}p_{i_g,t} \nonumber \\ 
	&- \sum_{\mathclap{u:(m,u) \in Z}}( q_{m,u,t}^{\text{in}}-q_{u,m,t}^{\text{out}})= 0:\hat{\lambda}^{\text{G}}_{m,t}, \enspace \forall m \in M, \thinspace \forall t \in T, \label{Gas_bal_coup}	
	\end{align}
\end{subequations}
where $\Theta^{\text{D}}=\{\Theta^{\text{ED}};\Theta^{\text{GD}}\}$ is the set of primal optimization variables. Objective function \eqref{Obj_DA_coup} determines the day-ahead operating cost of the electricity and natural gas systems. The system operating cost stems from the power production cost of thermal power plants (TPPs, i.e., non-gas) and the natural gas production cost. Note that the power production cost of GFPPs is not included, since it would result in double counting it. The cost of operating GFPPs is explicitly associated with the natural gas system cost through the balancing equation \eqref{Gas_bal_coup}. In this formulation, the fuel consumption is treated as a variable and charged with the locational marginal price ($\hat{\lambda}^{\text{G}}_{m,t}$) since the operation is jointly optimized. Having determined the day-ahead schedule of the integrated energy system, the real-time operation is simulated for each wind power realization $W_{j,\omega',t}$ by model \eqref{El_Gas_RT} that writes as follows,
\begin{subequations}
\label{El_Gas_RT}
	\begin{align}
	& \underset{\Theta^{\text{R}}}{\text{Min.}}\enspace \sum_{t \in T} \Big( \sum_{k \in K}(C_{k}^{+}g_{k,\omega',t}^{+}-C_{k}^{-}g_{k,\omega',t}^{-}) + \sum_{i_c \in I_c}(C_{i_c}^{+}p_{i_c,\omega',t}^{+}-C_{i_c}^{-}p_{i_c,\omega',t}^{-}) \nonumber\\
	& + \sum_{s \in S}(C_{s}^{+}g_{s,\omega',t}^{\text{out},r}-C_{s}^{-}g_{s,\omega',t}^{\text{in},r}) + \sum_{n \in N}C^{\text{sh,E}}l^{\text{sh,E}}_{n,\omega',t} 
	+ \sum_{m \in M} C^{\text{sh,G}}l^{\text{sh,G}}_{m,\omega',t} \Big) \label{Obj_RT_coup} 
	\intertext{subject to}
& \text{BA EL constraints} \enspace \eqref{RT_def}-\eqref{RefNode_RT},\label{RT_constraints_coup1}\\
& \text{BA NG constraints} \enspace \eqref{RT_capacityG1}-\eqref{Gas_shed},\eqref{RT_constraits_set1},\eqref{RT_constraits_set2},\label{RT_constraints_coup2}
\end{align}
\begin{align}
& \sum_{\mathclap{k \in A^{K}_{m}}} (g_{k,\omega',t}^{+} - g_{k,\omega',t}^{-})+ \sum_{\mathclap{s \in A_{m}^{S}}}( g_{s,\omega',t}^{\text{out,r}} - g_{s,\omega',t}^{\text{in,r}})  + \sum_{\mathclap{u:(m,u) \in Z}}(q_{m,u,t}^{\text{in,*}}-q_{u,m,t}^{\text{out,*}} - q_{m,u,t}^{\text{in,r}}+q_{u,m,t}^{\text{out,r}}) \nonumber \\
& - \sum_{\mathclap{i_g \in I_g}} \phi_{i_g}\Delta p_{i_g,\omega',t}  + l^{\text{sh,G}}_{m,\omega',t} = 0:\tilde{\lambda}^{\text{G}}_{m,\omega',t}, \enspace \forall m \in M, \thinspace \forall t \in T,  \label{RT_Gas_balance_coup}
\end{align}
\end{subequations}
where $\Theta^{\text{R}}=\{\Theta^{\text{ER}};\Theta^{\text{GR}}\}$ is the set of primal optimization variables. The goal is to minimize the cost of re-dispatch actions for the integrated energy system described in \eqref{Obj_RT_coup}. Similarly, the power production costs related to GFPPs are not included since they stem from the natural gas system. Fuel consumption of GFPPs is a variable in the co-optimization model of electricity and natural gas systems. Solving models \eqref{El_Gas_DA} and \eqref{El_Gas_RT} results in calculating the total expected cost under scenario set $\Omega$.

\subsection{Stochastic coupled electricity and natural gas model}
\label{StochCoup}
The stochastic and coupled dispatch model optimizes the operation of the integrated energy systems and is formulated in such a way that the day-ahead decisions anticipate future balancing costs over the scenario set $\Omega$. We formulate a two-stage stochastic programming problem \eqref{SC} that writes as follows,
\begin{subequations}
\label{SC}
	\begin{align}
	& \underset{\Theta^{\text{SC}}}{\text{Min.}}\enspace \sum_{t \in T} \Big[ \sum_{i_c \in I_c}C_{i_c}p_{i_c,t} + \sum_{k \in K}C_{k}g_{k,t} + \sum_{s \in S}C_{s}g_{s,t}^{\text{out}} \nonumber\\ \nonumber
	& + \sum_{\omega \in \Omega} \pi_{\omega} \Big( \sum_{k \in K}(C_{k}^{+}g_{k,\omega,t}^{+}-C_{k}^{-}g_{k,\omega,t}^{-}) + \sum_{i_c \in I_c}(C_{i_c}^{+}p_{i_c,\omega,t}^{+}-C_{i_c}^{-}p_{i_c,\omega,t}^{-}) \\ 
	&    + \sum_{s \in S}(C_{s}^{+}g_{s,\omega,t}^{\text{out},r}-C_{s}^{-}g_{s,\omega,t}^{\text{in},r}) + \sum_{n \in N}C^{\text{sh,E}}l^{\text{sh,E}}_{n,\omega,t} + \sum_{m \in M} C^{\text{sh,G}}l^{\text{sh,G}}_{m,\omega,t} \Big) \Big] \label{Obj_stoch}
\intertext{subject to} 
& \text{DA coupled EL \& NG constraints} \thinspace \eqref{El_max_min},\eqref{El_bal}\!-\!\eqref{RefNode_DA},\!\eqref{DA_constraints_coup2},\!\eqref{Gas_bal_coup}, \label{DA_STOCH}\\
& 0 \le w_{j,t} \le \overline{W}_j,\enspace \forall j \in J, \thinspace \forall t \in T, \label{wind_stoch} \\
& \text{BA coupled EL \& NG constraints} \enspace \eqref{RT_constraints_coup1}-\eqref{RT_Gas_balance_coup}, \enspace \forall \omega \in \Omega, \label{RT_STOCH}
\end{align}
\end{subequations}
where $\Theta^{\text{SC}}=\{\Theta^{\text{D}};\Theta^{\text{R}}_{\omega}, \thinspace \forall \omega \in \Omega\}$ is the set of primal optimization variables. The expected cost of operating the integrated energy system is given by the objective function \eqref{Obj_stoch} to be minimized. The stochastic and coupled dispatch model permits an implicit temporal coordination of the day-ahead and balancing stages through the expected balancing costs in \eqref{Obj_stoch} and constraints \eqref{RT_STOCH} that are modeled for all scenarios $\omega \in \Omega$. It should be noted that the day-ahead decisions are treated as variables and that the day-ahead wind power production is limited by the installed capacity through \eqref{wind_stoch}. This formulation optimally dispatches the integrated system to account for the uncertain wind power production by anticipating future balancing needs. Flexible producers may be scheduled out of merit-order to make flexibility available at the balancing stage. Moreover, the inclusion of natural gas system constraints in both day-ahead and balancing stages allows for an optimal spatial allocation of linepack flexibility in the system that is significantly important when getting closer to real-time operation.

\section{Results}
\label{Results}
\begin{figure}[ht]%
	\centering
	\includegraphics[scale=0.4]{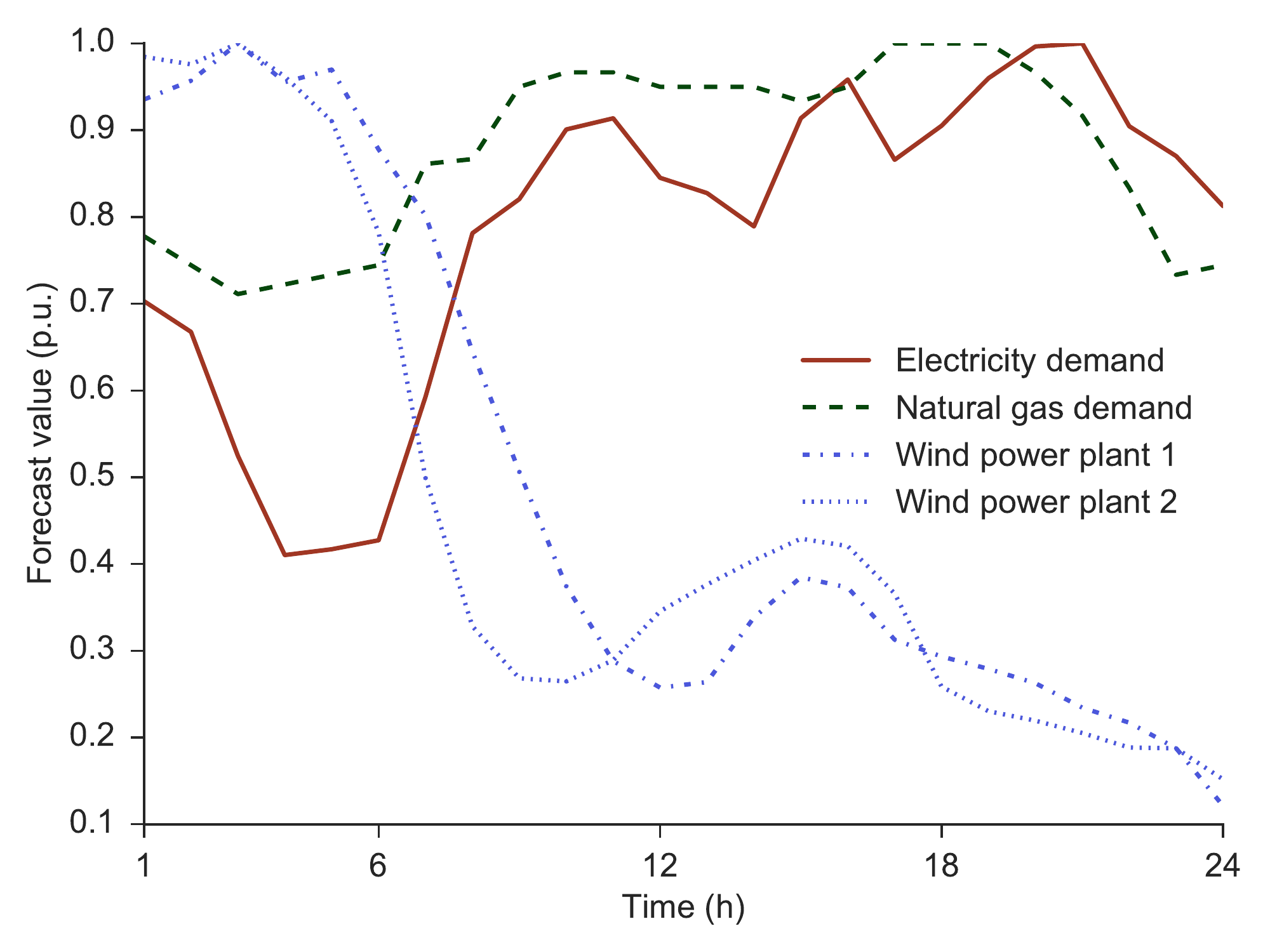}%
	\caption{Forecast profile of electricity and natural gas demands along with expected wind power production of wind power plant 1 and 2.}%
	\label{ForecastValuesFig}%
\end{figure}

A modified 24-bus IEEE Reliability Test System and a 12-node natural gas system compose the integrated energy system. It consists of 12 power plants, 2 wind farms, 17 electricity loads, 3 natural gas suppliers, 4 natural gas loads and 2 natural gas storage facilities. There are 4 flexible GFPPs that account for $427 \thinspace \text{MW}$ of the total $3075\thinspace \text{MW}$ installed capacity of conventional generation and interconnect the two systems. Wind power uncertainty is modeled by a set of 25 equiprobable scenarios that have proper temporal and spatial correlation, which are available at \citep{Scenarios}. The data of the integrated energy system is based on \citep{He2016} and presented along with the network topology in the online appendix available at \citep{Ordoudis2017a}. The forecast profile of electricity and residential natural gas demands, as well as the expected wind power production are illustrated in Fig. \ref{ForecastValuesFig}.

\subsection{Comparison of market-clearing models and the effect of coordination parameters}
In the following analysis, we highlight the inefficiencies arising from fixing the coordination parameters between the two markets in the decoupled approach and the benefits of dispatching the energy system in view of future uncertainties along with the linepack flexibility. We consider a 24-hour scheduling horizon, while the level of total system linepack and natural gas in the storage facilities at the end of day is equal to the one at the beginning and accounts for the operation of the following day. For the total system linepack, we set this value equal to $448\thinspace000\thinspace\text{kcf}$.

Initially, we examine the effect of utilizing fixed natural gas prices that are usually different from the natural gas marginal prices at each node of the system. Fixed natural gas prices represent the estimation of the natural gas spot price or the price of supply contracts that is used by GFPPs to bid in the electricity market. In our approach, the nodal natural gas prices stemming from the coupled clearing of electricity and natural gas markets are considered an appropriate estimation of the actual natural gas price. The solution of \textit{Seq-Coup} results in day-ahead natural gas prices ($\hat{\lambda}_{m,t}^{\text{G}}$) that are used in the \textit{Seq-Dec} model to define the power production cost of GFPPs. Specifically, the day-ahead offer price of GFPPs in \textit{Seq-Dec} is calculated by the multiplication of natural gas price of \textit{Seq-Coup} and the power conversion factor of each plant, i.e., $C_{i_g,t}=\hat{\lambda}_{m,t}^{\text{G}}\phi_{i_g}, \forall m\in A^{i_g}_m, \thinspace \forall i_g \in I_g, \thinspace \forall t \in T$. Upward and downward regulation offer prices are equal to 1.1 and 0.91 of the day-ahead offer price. A mis-estimation of the actual natural gas price by GFPPs is introduced in order to simulate the case that GFPPs bid in the electricity market with a natural gas price different that the actual one. This mis-estimation is simulated by a 10\% over- and under-estimation, while such deviation is considered adequate due to the relatively stable natural gas prices in short-term operations. The expected system costs are calculated based on \eqref{Obj_stoch} in order to ensure consistency of the results. Therefore, the offer prices of GFPPs in $\textit{Seq-Dec}$ only affect their position in the merit-order and the unit dispatch, while these prices are not taken into account in the cost calculation. 

Tables \ref{CostWind40} and \ref{CostWind50} show the costs of operating the energy system under the three market-clearing models and the share of day-ahead power production for GFPPs and TPPs for wind power penetration levels, i.e., share of installed wind power capacity on system's demand, of 40\% and 50\%, respectively. It is observed that \textit{Stoch-Coup} achieves the lowest expected cost in both cases by efficiently accommodating the large shares of renewable power production. This model decides the optimal day-ahead dispatch to account for wind power uncertainty which results in a higher day-ahead cost but lower expected balancing cost, while the share of GFPPs that are efficient balancing producers also increases compared to deterministic models \textit{Seq-Coup} and \textit{Seq-Dec}. The deterministic models schedule the system based on the merit-order principle which may not be appropriate to cope with wind power uncertainty. We would like to note that the characteristics of wind power scenarios (e.g., mean and variance) affect the outcome of the market-clearing models. The \textit{Stoch-Coup} model handles uncertainty more efficiently than the deterministic ones and this can be observed by the greater decrease in expected cost it accomplishes when increasing the wind power penetration level from 40\% to 50\%. The initial wind power scenarios are normalized and then multiplied by the wind farm capacity; thus a higher penetration level results in scenarios with higher mean and variance. We refer the reader to \citep{Ordoudis2016} for further discussion on this issue. 

A coupled clearing of electricity and natural gas markets results in lower expected cost compared to the decoupled approaches due to the optimal natural gas price signals, as well as the optimized natural gas flows and linepack. The spatial allocation of natural gas in the system plays an important role for short-term adequacy and available flexibility in view of the ability to store gas in the pipelines. For the cases of natural gas price mis-estimation, it can be noticed that there might be different effects on the total expected cost depending on wind power penetration level. The day-ahead cost increases when the natural gas price is mis-estimated compared to \textit{Seq-Dec} that utilizes the actual natural gas price from \textit{Seq-Coup}. The reason for that is the misplacement of GFPPs in the merit-order because of the over-estimation ($\uparrow$) or under-estimation ($\downarrow$) of natural gas price compared to the merit-order built by the actual one. This results in scheduling more expensive units at the day-ahead stage  which in turn increases the day-ahead cost. In the under-estimation cases, flexible GFPPs are scheduled more at the day-ahead stage which makes them unavailable to provide up-regulation in real-time operation, thus the expected balancing cost significantly increases as more expensive sources have to be deployed. On the contrary, it is possible to have a reduced total expected cost when the natural gas price is over-estimated. In the case of 50\% wind power penetration, flexible GFPPs are scheduled less at the day-ahead stage which makes them available to provide upward regulation services in real-time.  

\begin{table}[ht]
\centering
\caption{Expected cost and share of day-ahead power production. Wind power penetration 40\%.}
\label{CostWind40}
\resizebox{0.75\columnwidth}{!}{
\begin{tabular}{cccccc}
\hline
           & Total (\$)                  & Day-ahead (\$)              & Balancing (\$)   & GFPPs (\%) & TPPs (\%) \\ \hline
Stoch-Coup & 1\thinspace747\thinspace156 & 1\thinspace755\thinspace474 & -8\thinspace318  & 10.8             & 89.2            \\
Seq-Coup   & 1\thinspace817\thinspace781 & 1\thinspace731\thinspace522 & 86\thinspace259  & 9.1              & 90.9            \\
Seq-Dec & 1\thinspace819\thinspace238 & 1\thinspace731\thinspace657 & 87\thinspace581  & 9.2              & 90.8            \\
Seq-Dec $\uparrow$ & 1\thinspace821\thinspace399 & 1\thinspace732\thinspace185 & 89\thinspace214  & 6.8              & 93.2            \\
Seq-Dec $\downarrow$ & 1\thinspace866\thinspace966 & 1\thinspace741\thinspace471 & 125\thinspace495 & 14.5             & 85.5            \\ \hline
\end{tabular}
}
\end{table}

\begin{table}[ht]
\centering
\caption{Expected cost and share of day-ahead power production. Wind power penetration 50\%.}
\label{CostWind50}
\resizebox{0.75\columnwidth}{!}{
\begin{tabular}{cccccc}
\hline
                        & Total (\$)                  & Day-ahead (\$)              & Balancing (\$)   & GFPPs (\%) & TPPs (\%) \\ \hline
Stoch-Coup              & 1\thinspace684\thinspace075 & 1\thinspace686\thinspace504 & -2\thinspace429  & 11.7       & 88.3      \\
Seq-Coup                & 1\thinspace812\thinspace790 & 1\thinspace663\thinspace477 & 149\thinspace313 & 7.4        & 92.6      \\
Seq-Dec              & 1\thinspace814\thinspace086 & 1\thinspace663\thinspace891 & 150\thinspace195 & 7.2        & 92.8      \\
Seq-Dec $\uparrow$   & 1\thinspace813\thinspace271 & 1\thinspace664\thinspace646 & 148\thinspace625 & 6.1        & 93.9      \\
Seq-Dec $\downarrow$ & 1\thinspace874\thinspace907 & 1\thinspace674\thinspace110 & 200\thinspace797 & 14.3       & 85.7      \\ \hline
\end{tabular}
}
\end{table}

\begin{figure}[ht]%
\vspace{-7.5pt}
	\centering
	\includegraphics[scale=0.4]{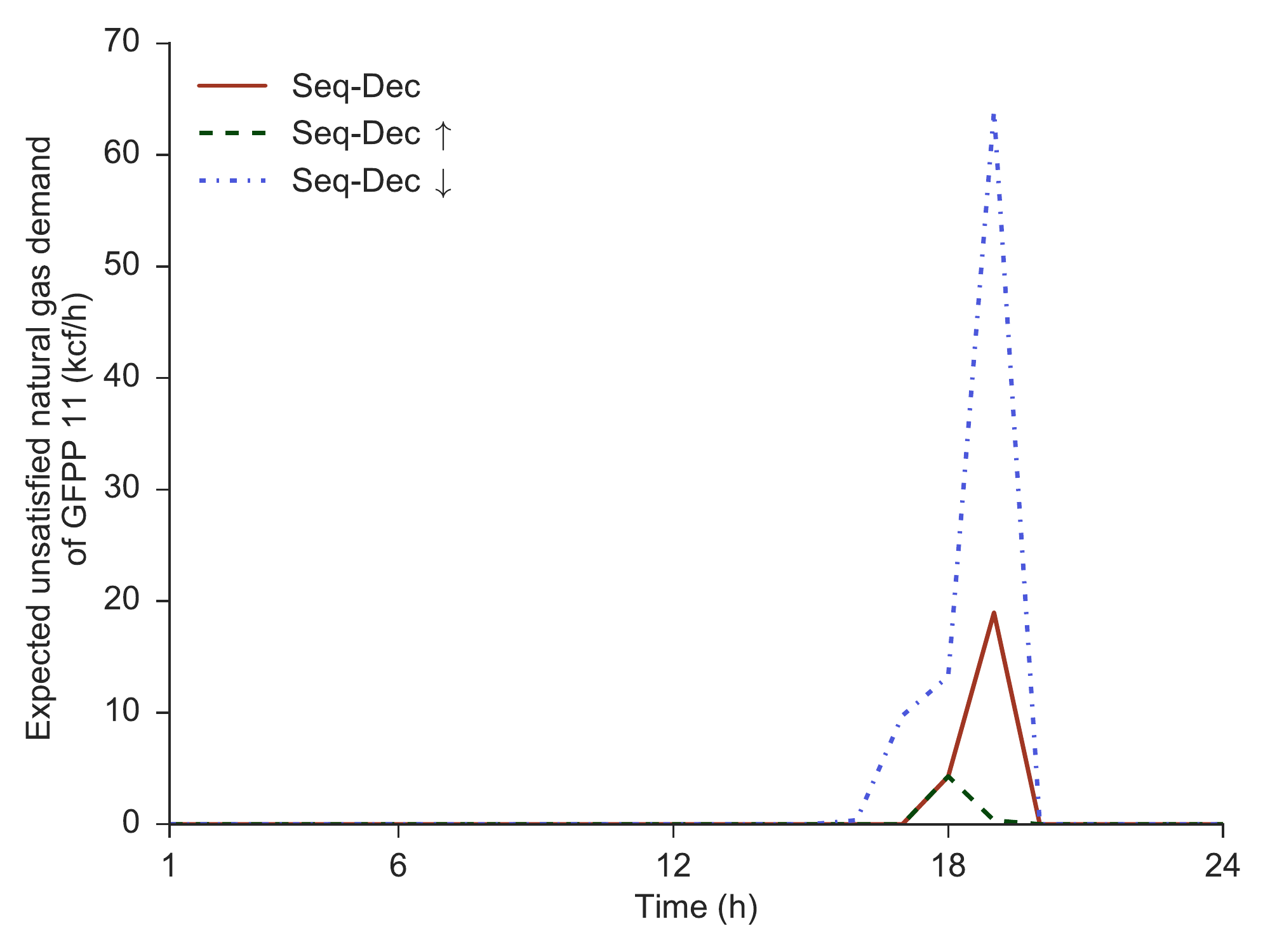}%
	\caption{Expected unsatisfied natural gas demand of GFPP 11 under Seq-Dec.}%
	\label{UnsatisfiedDemandFig}%
\end{figure}

Additionally, the impact of congestion in the natural gas network on the scheduling of GFPPs is illustrated in Table \ref{CostNGDWind40} and Fig. \ref{UnsatisfiedDemandFig}. The residential natural gas demand is increased by 30\% to represent winter conditions, when natural gas is used extensively for heating. The share of GFPPs is reduced from 9.1\% to 6.6\% and from 9.2\% to 6.7\% for \textit{Seq-Coup} and $\textit{Seq-Dec}$, respectively. The TPPs have to cover a higher portion of the electricity demand in this case due to the lower priority assigned to the fuel demand of GFPPs compared to  natural gas residential loads. Moreover, it is observed that GFPP 11 is expected to face unsatisfied fuel demand during the hours of peak residential natural gas load under $\textit{Seq-Dec}$. This phenomenon does not occur when the systems are simultaneously operated in $\textit{Seq-Coup}$. In this case, the importance of efficiently optimizing the gas flows and the spatial allocation of linepack in the system is more evident. The expected unsatisfied natural gas demand of GFPP 11 is higher when the natural gas price is under-estimated as it is scheduled more at the day-ahead stage and lower when the natural gas price is over-estimated. 

\begin{table}[ht]
\centering
\caption{Expected cost and share of day-ahead power production. Wind power penetration 40\%. Increased natural gas demand by 30\%.}
\label{CostNGDWind40}
\begin{tabular}{cccccc}
\hline
           & Total (\$)                  & GFPPs (\%) & TPPs (\%) \\ \hline
Seq-Coup   & 1\thinspace949\thinspace993  & 6.6        & 93.4      \\ 
Seq-Dec   & 2\thinspace018\thinspace479 & 6.7        & 93.3      \\ \hline
\end{tabular}
\end{table}

\subsection{The benefits and effects of linepack}
The following results aim at exploring the benefits of linepack and the corresponding flexibility revealed. In the following analysis, the wind power penetration level is 50\%. Initially, we simulate a purely steady-state operation of the natural gas system, where no linepack is considered and the inflow and outflow of each pipeline is equal ($q_{m,u}^{\text{in}}=q_{m,u}^{\text{out}}, \forall (m,u) \in Z$).

Table \ref{LinepackCostSeq} presents the expected system cost under \textit{Seq-Coup} and \textit{Seq-Dec} comparing the cases where linepack is considered or neglected (steady-state). The steady-state models result in a higher day-ahead cost for both \textit{Seq-Coup} and \textit{Seq-Dec} because the most expensive natural gas producer is scheduled to meet the demand. In this case, the natural gas demand has to be instantly covered at each time period by the production units as it is not possible to store natural gas in the network for future utilization. When the linepack is considered, natural gas is only supplied by the two cheaper producers and it is stored in the beginning of the day to be used later during the hours of peak demand. It is observed, though, that the expected balancing costs are lower in the steady-state models. The deployment of the expensive natural gas supplier increases the natural gas price and makes available more cost-effective capacity for down-regulation in the natural gas system. This increased availability of cheaper down-regulation is also reflected in the electricity market through the GFPPs and results in lower expected balancing cost. This observation depends on the structure of scenarios and the type of regulation needed in real-time operation. Nevertheless, it indicates possible inefficiencies that may arise when a flexible component of the system, such as the linepack, is myopically operated. 
\begin{table}[ht]
\centering
\caption{Effect of linepack on expected cost under \textit{Seq-Coup} and \textit{Seq-Dec}. Wind power penetration 50\%.}
\label{LinepackCostSeq}
\begin{tabular}{cccc}
\hline
                  & Total (\$)                  & Day-ahead (\$)              & Balancing (\$)   \\ \hline
Seq-Coup          & 1\thinspace812\thinspace790 & 1\thinspace663\thinspace477 & 149\thinspace313 \\
Steady-state Seq-Coup   & 1\thinspace809\thinspace753 & 1\thinspace669\thinspace452 & 140\thinspace301 \\
Seq-Dec        & 1\thinspace814\thinspace086 & 1\thinspace663\thinspace891 & 150\thinspace195 \\
Steady-state Seq-Dec & 1\thinspace813\thinspace271 & 1\thinspace671\thinspace109 & 142\thinspace162 \\ \hline
\end{tabular}
\end{table}

The subsequent results show that \textit{Stoch-Coup} effectively schedules the system and exploits the available flexibility of linepack under uncertainty of wind power production. The temporal coordination between the two trading floors allows an efficient allocation of natural gas in the system and overcome the aforementioned drawback of the deterministic approaches that are myopic to uncertainties. The steady-state model results in higher expected cost compared to $\textit{Stoch-Coup}$ when linepack is taken into account, as illustrated in Table \ref{LinepackCostStoch}. Moreover, $\textit{Stoch-Coup}$ schedules more the GFPPs by taking advantage of the storage ability in the natural gas system. 
\begin{table}[ht]
\centering
\caption{Effect of linepack on expected cost and share of day-ahead power production under \textit{Stoch-Coup}. Wind power penetration 50\%.}
\label{LinepackCostStoch}
\resizebox{0.75\columnwidth}{!}{
\begin{tabular}{cccccc}
\hline
                  & Total (\$)                  & GFPPs (\%) & TPPs (\%) \\ \hline
Steady-state Stoch-Coup  & 1\thinspace691\thinspace728 & 9.9        & 90.1      \\
Stoch-Coup        & 1\thinspace684\thinspace076 & 11.7       & 88.3      \\
Stoch-Coup (+5\% initial linepack) & 1\thinspace631\thinspace559 & 12         & 88        \\
Stoch-Coup (-5\% initial linepack) & 1\thinspace739\thinspace304 & 10.8       & 89.2      \\ \hline
\end{tabular}
}
\end{table}

The main effect of considering a storage facility in the operation of the energy system is that it flattens the demand profile by filling valleys and shaving peaks, as well as utilizing cheap power production. In order to quantify the flexibility revealed by modeling linepack in the integrated energy system, we simulate a case where an ideal storage facility, i.e., having infinite capacity and charging/discharging rates, is introduced in the electricity network allowing to shift the demand profile in the most cost-effective manner. In this case, the resulting expected cost ($EC$) of $\textit{Stoch-Coup}$ is $\$1\thinspace629\thinspace519$. The linepack flexibility is quantified by the following performance ratio,
\begin{equation}
\begin{aligned}
    &\frac{EC^{ss}-EC}{EC^{ss}-EC^{ideal}}=\frac{1\thinspace691\thinspace728-1\thinspace684\thinspace016}{1\thinspace691\thinspace728-1\thinspace629\thinspace519}=12.4 \%. \nonumber
\end{aligned}
\end{equation}
This ratio shows that modeling the natural gas network with linepack flexibility unveils 12.4\% of utilizing an ideal storage in the electricity system. This result indicates the need of efficiently modeling the integrated energy market as it provides an economical approach to exploit the available flexibility in the natural gas system. 

Furthermore, we examine the effect of the linepack level at the beginning of the day on the system operation by defining two additional cases for the ease of our analysis. These two cases are defined by having a 5\% more or less linepack at the beginning of the day in relation to the value at the end of the day. In all cases, the total linepack at the end of the day is equal to $448\thinspace000\thinspace\text{kcf}$. The total operating cost is lower when there is higher level of initial linepack in the system. The system is operated by taking advantage of the free energy stored in the natural gas network that also results in scheduling more the GFPPs compared to the base case of \textit{Stoch-Coup}. On the contrary, the expected cost is higher when the gas network has to be filled in order to meet the final condition of the scheduling horizon.

\begin{figure}[ht]%
\vspace{-7.5pt}
	\centering
	\includegraphics[scale=0.4]{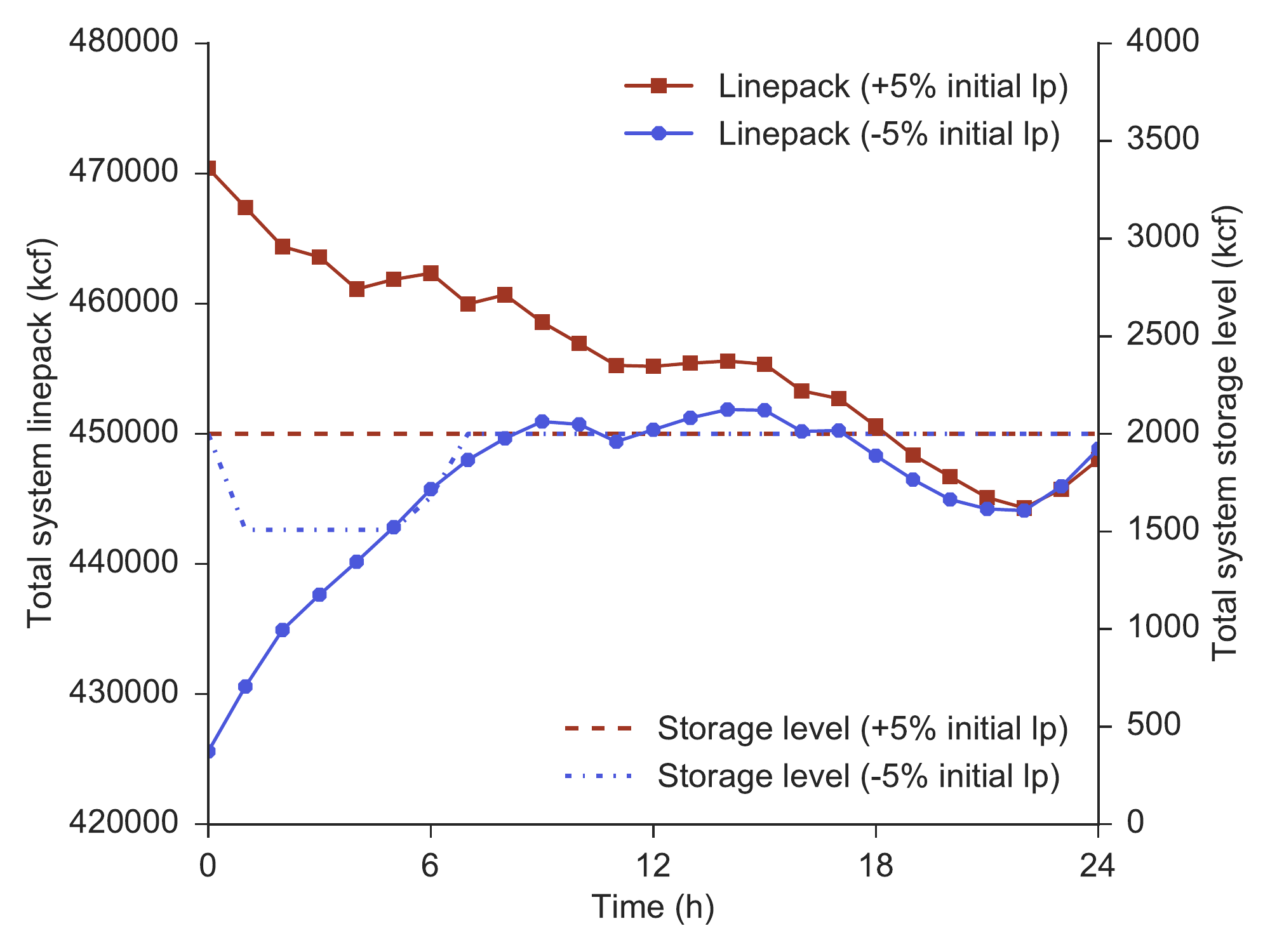}%
	\caption{Total system linepack and storage level (lp: linepack).}%
	\label{LinepackStorageFig}%
\end{figure}

The total linepack of the system is presented in Fig. \ref{LinepackStorageFig}, along with the total storage level. The linepack decreases throughout the scheduling horizon when the initial level is higher than the final condition and the storage is not utilized at the day-ahead stage. It can be noticed that the rate of decrease depends on the residential natural gas demand profile, while there are couple of periods that the linepack is charged when the residential demand is relatively low. On the other hand, the linepack is increased when its initial value is lower at the beginning of the scheduling horizon. The highest rate of increase is observed when the residential natural gas demand is low during the first hours of the day. Moreover, there is a need to discharge the storage facilities to support the secure operation of the system. In both cases, the linepack is decreased below the final hour threshold during the hours of peak residential demand to avoid supplying natural gas by the expensive natural gas producer.

\begin{figure}[ht]%
\vspace{-7.5pt}
	\centering
	\includegraphics[scale=0.4]{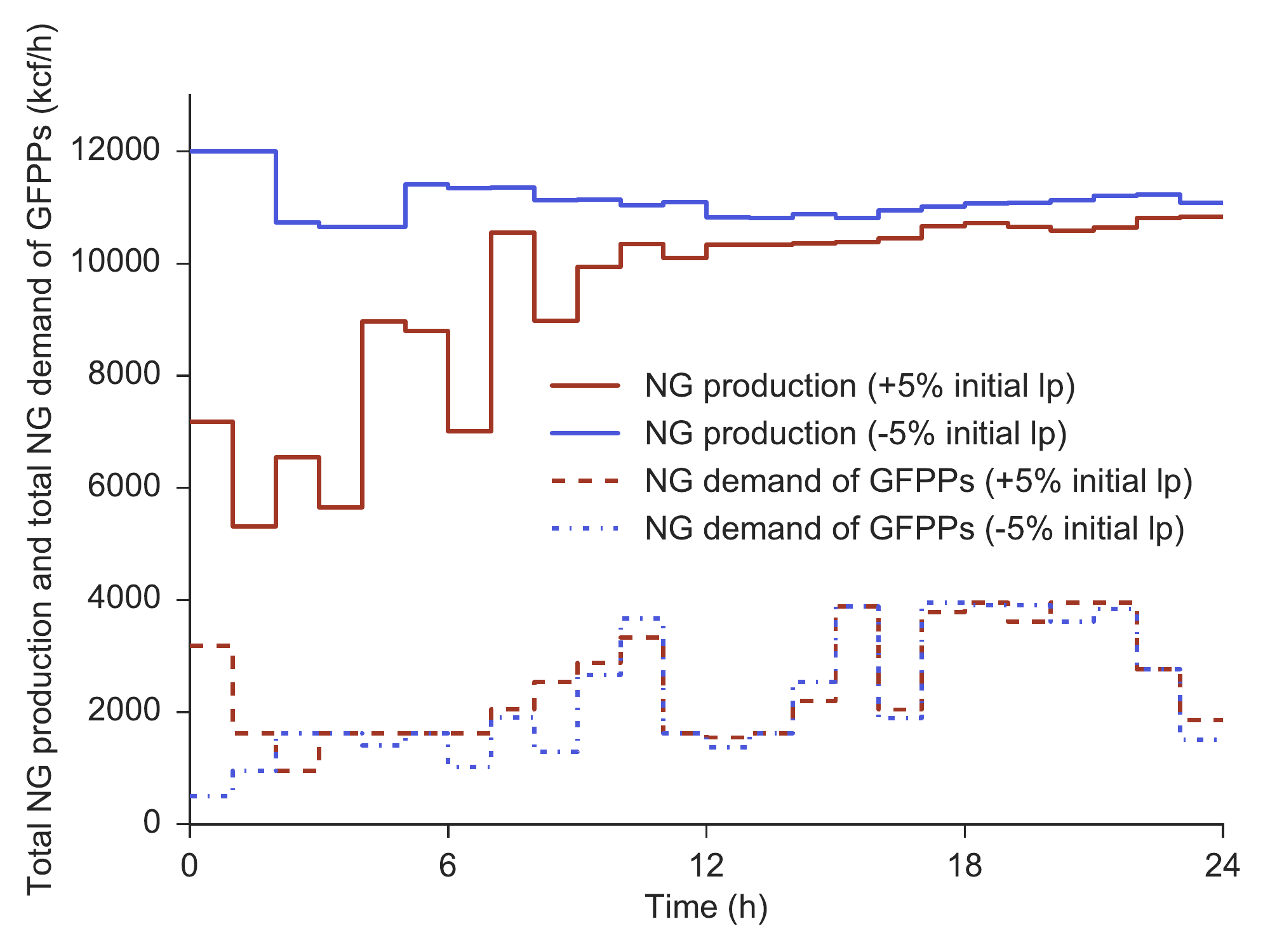}%
	\caption{Total natural gas production and total natural demand of GFPPs (NG: natural gas, lp: linepack).}%
	\label{GasProductionFigGFPPsDemandFig}%
\end{figure}

In addition, Fig. \ref{GasProductionFigGFPPsDemandFig} illustrates the total natural gas supply and the fuel demand of GFPPs. The natural gas production is significantly higher when the lack of linepack in the system has to be covered. On the contrary, the utilization of linepack to cover natural gas demand during the first half of the scheduling horizon, when the initial value is higher, is demonstrated by the reduced production in the corresponding hours. Moreover, the GFPPs are scheduled more in this case due to the excess of natural gas in the system. This points out another advantage of the coupled approach where the GFPPs will be scheduled according to the gas availability in the system. The GFPPs are serving as a flexible demand component for the natural gas system, by either increasing or decreasing their fuel consumption, in favor of a cost-effective operation. 

The optimization problems were solved using CPLEX 12.6.2 under GAMS on a stationary computer with Intel i7 4-core processor clocking at 3.4 GHz and 8 GB of RAM. The average time to solve \textit{Stoch-Coup} was 1540 seconds, while \textit{Seq-Coup} and \textit{Seq-Dec} were solved in less than 90 seconds. 

\section{Conclusion}
\label{Conclusion}

This paper proposes a co-optimization model for integrated electricity and natural gas systems that efficiently takes into account uncertain power supply. We follow a linearization approach to approximate the dynamics of the natural gas system that yields a tractable mixed-integer linear programming (MILP) model and considers the possibility to store gas in the pipelines of the natural gas network (i.e., linepack), which is of significant importance in short-term operations. The combination of the aforementioned model properties results in increasing operational flexibility and in improved allocation of gas resources in the network. Moreover, we provide two market-clearing models with a deterministic description of the uncertain power supply and assess the value of coordination between the energy systems and the trading floors. 

Initially, the impact of coordination parameters, such as consumption of gas-fired power plants (GFPPs) and natural gas price, is examined in a setup where the electricity and natural gas markets are decoupled. It is shown that both parameters affect the operation of both systems and in turn, the expected cost. A poor definition of these parameters between these markets usually has a negative effect on the expected cost. However, a distorted day-ahead scheduling of GFPPs due to arbitrary definition of these coordination parameters may result in making available reserve capacity that is suitable to cope with the actual imbalances in real-time operation. In contrast, a fully coupled model effectively copes with power supply uncertainty and makes the most of the flexibility inherent to the natural gas system, which is further highlighted by comparing the approaches that linepack is modeled or disregarded. Moreover, our analysis shows that it is not adequate to introduce an additional source of flexibility (e.g., linepack) in the system without having the proper market design to operate it. Finally, simulation results show that the models efficiently capture the natural gas system behavior.

For future works, the optimal definition of coordination parameters between the two markets needs to be studied. Results in our analysis indicate that a systematic approach to define them is capable of reducing the expected system cost. Moreover, a more detailed description of the natural gas system that models also fuel consumption of compressors could be considered. This would highly increase the computational complexity of the problem. However, tailored decomposition techniques can be applied to overcome this hurdle.  
\section*{Acknowledgment}
The work of Christos Ordoudis and Pierre Pinson is partly funded by the Danish Strategic Research Council (DSF) through projects “5s-Future Electricity Markets”, No. 12-132636/DSF and “CITIES”, No. 1305-00027B/DSF. The work of Juan M. Morales is partly funded by the European Research Council (ERC) under the European Union’s Horizon 2020 research and innovation programme (grant agreement No. 755705), by the Spanish Research Agency through project ENE2017-83775-P (AEI/FEDER, UE), and by the Research Funding Program for Young Talented Researchers of the University of M\'{a}laga through project PPIT-UMA-B1-2017/18. Finally, we would like to thank the editor and the reviewers for their comments to improve this manuscript.

     \printnomenclature

\nomenclature[S]{$I$}{Set of dispatchable power plants $i$.}
\nomenclature[S]{$I_c$}{Set of thermal power plants $i_c$ ($I_c \subset I$).}
\nomenclature[S]{$I_g$}{Set of natural gas-fired power plants $i_g$ ($I_g \subset I$).}
\nomenclature[S]{$J$}{Set of stochastic power plants $j$.}
\nomenclature[S]{$L$}{Set of electricity transmission lines $l$.}
\nomenclature[S]{$N$}{Set of electricity network nodes $n$.}
\nomenclature[S]{$K$}{Set of natural gas producers $k$.}
\nomenclature[S]{$S$}{Set of natural gas storages $s$.}
\nomenclature[S]{$Z$}{Set of natural gas network branches $z$.}
\nomenclature[S]{$M$}{Set of natural gas network nodes $m$.}
\nomenclature[S]{$V$}{Set of fixed pressure points $v$ for the linearization of Weymouth equation.}
\nomenclature[S]{$\Omega$}{Set of stochastic power production scenarios $\omega$.}
\nomenclature[S]{$T$}{Set of time periods $t$.}
\nomenclature[S]{$A_{n}^{I}$}{Set of dispatchable power plants $i$ located at electricity network node $n$.}
\nomenclature[S]{$A_{n}^{J}$}{Set of stochastic power plants $j$ located at electricity network node $n$.}
\nomenclature[S]{$A_{m}^{I_g}$}{Set of natural gas-fired power plants $i_g$ located at natural gas network node $m$.}
\nomenclature[S]{$A_{m}^{K}$}{Set of natural gas producers $k$ located at natural gas network node $m$.}
\nomenclature[S]{$A_{m}^{S}$}{Set of natural gas storages $s$ located at natural gas network node $m$.}
\nomenclature[S]{$\Theta$}{Set of primal optimization variables defined for each optimization model.}
\nomenclature[P]{$D^{\text{E}}_{n,t}$}{Electricity demand at node $n$ and in period $t$ [MW]. } 
\nomenclature[P]{$D^{\text{G}}_{m,t}$}{Natural gas demand at node $m$ and in period $t$ [kcf/h].} 
\nomenclature[P]{$C_{i}$}{Day-ahead offer price of dispatchable power plant $i$ [\$/MWh].}
\nomenclature[P]{$C^{+}_{i},C^{-}_{i}$}{Up/down regulation offer price of dispatchable power plant $i$ [\$/MWh].}
\nomenclature[P]{$C^{\text{sh,E}}$}{Cost of electricity load shedding [\$/MWh].}
\nomenclature[P]{$C_{k}, C_{s}$}{Day-ahead offer price of natural gas producer $k$ and storage $s$ [\$/kcf].}
\nomenclature[P]{$C^{+}_{k},C^{-}_{k}$}{Up/down regulation offer price of natural gas producer $k$ [\$/kcf].}
\nomenclature[P]{$C^{+}_{s},C^{-}_{s}$}{Up/down regulation offer price of natural gas storage $s$ [\$/kcf].}
\nomenclature[P]{$C^{\text{sh,G}}$}{Cost of natural gas load shedding [\$/kcf].}
\nomenclature[P]{$P^{\text{max}}_i$}{Capacity of dispatchable power plant $i$ [MW].}
\nomenclature[P]{$P^+_i,P^-_i$}{ Maximum up/down reserve offered by dispatchable power plant $i$ [MW].}
\nomenclature[P]{$\phi_{i_g}$}{Power conversion factor of natural gas-fired power plant $i_g$ [kcf/MWh].}
\nomenclature[P]{$W_{j,\omega,t}$}{Power production by stochastic power plant $j$ in scenario $\omega$, period $t$ [MW].}
\nomenclature[P]{$\widehat{W}_{j,t}$}{Expected power production by stochastic power plant $j$ in period $t$ [MW].}
\nomenclature[P]{$\overline{W}_{j}$}{Capacity of stochastic power plant $j$ [MW].}
\nomenclature[P]{$G_{k}^{\text{max}}$}{Capacity of natural gas producer $k$ [$\text{kcf/h}$].}
\nomenclature[P]{$G_{k}^+,G_{k}^-$}{Maximum up/down reserve offered by natural gas producer $k$ [$\text{kcf/h}$].}
\nomenclature[P]{$B_{n,r}$}{Absolute value of the susceptance of line (n,r) [per unit].}
\nomenclature[P]{$F^{\text{max}}_{n,r}$}{Transmission capacity of line (n,r) [MW].}  
\nomenclature[P]{$E_{s}^{\text{min/max}}$}{Minimum and maximum level of storage $s$ [$\text{kcf}$].}
\nomenclature[P]{$K^{\text{h/f}}_{m,u}$}{Linepack (h) and natural gas flow (f) constant of pipeline (m,u) [kcf/psig, kcf/(psig $\cdot$ h)].}
\nomenclature[P]{$IR_{s},WR_{s}$}{Injection and withdrawal rates of storage $s$ [$\text{kcf/h}$].}
\nomenclature[P]{$PR_{m}^{\text{min/max}}$}{Minimum and maximum pressure at node $m$ [psig].}
\nomenclature[P]{$\Gamma_z$}{Compressor factor located at natural gas network branch $z$ [-].}
\nomenclature[P]{$\tilde{M}$}{Sufficiently large constant [-].}
\nomenclature[P]{$\pi_{\omega}$}{Probability of scenario $\omega$ [-].}
\nomenclature[V]{$p_{i,t},w_{j,t}$}{Day-ahead dispatch of power plants $i$ and $j$ in period $t$ [MW].}
\nomenclature[V]{$p^{+/-}_{i,\omega,t}$}{Up/down regulation provided by dispatchable power plant $i$ in scenario $\omega$, period $t$ [MW].}
\nomenclature[V]{$w^{\text{sp}}_{j,\omega,t}$}{Power spilled by stochastic power plant $j$ in scenario $\omega$, period $t$ [MW].}
\nomenclature[V]{$l^{\text{sh,E/G}}_{n/m,\omega,t}$}{Electric power and natural gas load shedding at node $n$/$m$ in scenario $\omega$, period $t$ [MW, kcf/h].}
\nomenclature[V]{$\hat{\delta}_{n,t}$}{Voltage angle at node $n$ and in period $t$ [rad].}
\nomenclature[V]{$\tilde{\delta}_{n,\omega,t}$}{Voltage angle at node $n$ in scenario $\omega$, period $t$ [rad].}
\nomenclature[V]{$g_{k,t}$}{Day-ahead dispatch of natural gas producer $k$ in period $t$ [kcf/h].}
\nomenclature[V]{$g^{+/-}_{k,\omega,t}$}{Up/down regulation provided by natural gas producer $k$ in scenario $\omega$, period $t$ [kcf/h].}
\nomenclature[V]{$pr_{m,t}$}{Pressure at node $m$ and in period $t$ [psig].}
\nomenclature[V]{$h_{m,u,t}$}{Average mass of natural gas (linepack) in pipeline (m,u), period $t$ [kcf].}
\nomenclature[V]{$g^{\text{in/out}}_{s,t}$}{In- and outflow natural gas rates of storage $s$ in period $t$ [kcf/h].}
\nomenclature[V]{$q^{\text{in/out}}_{m,u,t}$}{In- and outflow natural gas rates of pipeline (m,u) in period $t$ [kcf/h].}
\nomenclature[V]{$pr^{\text{r}}_{m,\omega,t}$}{Pressure at node $m$ in scenario $\omega$, period $t$ [psig].}
\nomenclature[V]{$h^{\text{r}}_{m,u,\omega,t}$}{Average mass of natural gas (linepack) in pipeline (m,u), scenario $\omega$, period $t$ [kcf].}
\nomenclature[V]{$g^{\text{in/out,r}}_{s,\omega,t}$}{In- and outflow natural gas rates of storage $s$ in scenario $\omega$, period $t$ [kcf/h].}
\nomenclature[V]{$q^{\text{in/out,r}}_{m,u,\omega,t}$}{In- and outflow natural gas rates of pipeline (m,u) in scenario $\omega$, period $t$ [kcf/h].}
\nomenclature[V]{$q_{m,u}$}{Natural gas flow in pipeline (m,u) [kcf/h].}
\nomenclature[V]{$y_{m,u,t}$}{Binary variable defining the direction of the natural gas flow in pipeline (m,u), period $t$ \{0,1\}.}
\nomenclature[V]{$y^{\text{r}}_{m,u,\omega,t}$}{Binary variable defining the direction of the natural gas flow in pipeline (m,u), scenario $\omega$, period $t$ \{0,1\}.}
\nomenclature[V]{$e_{s,t}$}{Natural gas volume in storage facility $s$ and in period $t$ [kcf].}
\nomenclature[V]{$\hat{f}_{n,r,t}$}{Power flow on line (n,r) and in period $t$ [MW].}
\nomenclature[V]{$\tilde{f}_{n,r,\omega,t}$}{Power flow on line (n,r), in scenario $\omega$, period $t$ [MW].}
\nomenclature[V]{$\hat{\lambda}^{\text{G}}_{m,t}$}{Natural gas locational marginal price in day-ahead market at node $m$ and period $t$ [\$/kcf].}
\nomenclature[V]{$\tilde{\lambda}^{\text{G}}_{m,\omega,t}$}{Natural gas locational marginal price in balancing market at node $m$ in scenario $\omega$, period $t$ [\$/kcf].}
\nomenclature[V]{$\hat{\lambda}^{\text{E}}_{n,t}$}{Electricity locational marginal price in day-ahead market at node $n$ and period $t$ [\$/MWh].}
\nomenclature[V]{$\tilde{\lambda}^{\text{E}}_{n,\omega,t}$}{Electricity locational marginal price in balancing market at node $n$ in scenario $\omega$, period $t$ [\$/MWh].}




\bibliographystyle{model5-names}\biboptions{authoryear}

\bibliography{Integrated_El_NG_new2}

\end{document}